\documentstyle[12pt,dina4,euler]{article}
\input amssym.def
\input amssym
\newcommand{\N}{{\Bbb N}}
\newcommand{\Z}{{\Bbb Z}}
\newcommand{\C}{{\Bbb C}}
\newcommand{\Q}{{\Bbb Q}}
\renewcommand{\P}{{\Bbb P}}
\newcommand{\kl}{{\cal L}}
\newcommand{\ko}{{\cal O}}
\newcommand{\ks}{{\cal S}}
\newcommand{\kx}{{\cal X}}

\begin{document}
\newtheorem{lemma}{Lemma}[section]
\newtheorem{remark}[lemma]{Remark}
\newtheorem{proposition}[lemma]{Proposition}
\newtheorem{example}[lemma]{Example}
\newtheorem{theorem}[lemma]{Theorem}
\newtheorem{corollary}[lemma]{Corollary}
\newtheorem{definition}[lemma]{Definition}

\title{{\bf Moduli spaces of semiquasihomogeneous singularities
with fixed principal part}}
\author{
Gert-Martin Greuel\\
Universit\"at Kaiserslautern\\
Fachbereich Mathematik\\
Erwin-Schr\"odinger-Stra\ss e\\
D -- 67663 Kaiserslautern
\and
Claus Hertling\\
Universit\"at Bonn\\
Mathematisches Institut\\
Wegelerstra\ss e 10\\
D -- 53115 Bonn
\and
Gerhard Pfister\\
Universit\"at Kaiserslautern\\
Fachbereich Mathematik\\
Erwin-Schr\"odinger-Stra\ss e\\
D -- 67663 Kaiserslautern}

\date{January 11,1995}
\maketitle
\vspace{1cm}
\tableofcontents
\thispagestyle{empty}
\setcounter{page}{0}
\newpage

\section*{Introduction}\addcontentsline{toc}{section}{Introduction}

One of the important achievements of singularity theory is the explicit
classification of certain ``generic'' classes of isolated hypersurface
singularities via normal forms and the analysis of its properties (cf.\
\cite{AGV}).   More complicated singularities deform into a
collection of singularities from these classes and deformation theory is a
powerful tool in studying specific singularities.   For a further
classification of more complicated classes of singularities the explicit
determination of normal forms seems to be impossible and not appropriate.

The aim of this article is to start towards a classification of isolated
hypersurface singularities of any dimension via geometric methods, that is by
explicitely constructing a (coarse) moduli space for such singularities with
certain
invariants being fixed.   Our method starts from deformation theory and leads
to the
construction of geometric quotients of quasiaffine spaces by certain algebraic
groups whose main part is unipotent.   This last part is a major ingredient
and uses the general results of \cite{GP 2}.   In
projective algebraic geometry, the theory of moduli spaces is highly developed
but in singularity theory only a few attempts have been made so far, for
example by Ebey, Zariski, Laudal, Pfister, Luengo, Greuel (cf.\ \cite{LP} for
a systematic approach and \cite{GP 1}
for a short survey).   In this paper we consider only
semiquasihomogeneous singularities given as a power series $f \in \C \{x_1,
\ldots, x_n\}$ or as a complex space germ $(X,0) = (f^{-1}(0),0) \subset
(\C^n,0)$, together with positive weights $w_1, \ldots, w_n$ of the variables
such that the principal part $f_0$ of $f$ (terms of lowest degree) has an
isolated singularity.

For the classification we first fix the Milnor number, probably the most basic
invariant of an isolated hypersurface singularity.   Fixing the Milnor number
is known (for $n \not= 3)$ to be equivalent
to fixing, in a family, the embedded topological type of the singularity.   If
the Milnor number is fixed, the classification of semiquasihomogeneous
singularities falls naturally into two parts.
Firstly, the classification of the quasihomogeneous principal parts or, which
amounts to the same, the classification of hypersurfaces in a weighted
projective space.   Secondly, the classification of semiquasihomogeneous
hypersurface singularities with fixed principal part.   These two parts differ
substantially, since the group actions whose orbits describe isomorphism
classes of singularities are of a completely different nature.   This article
is devoted to the second task.

The most important equivalence relations for hypersurface singularities are
right equivalence (change of coordinates in the source) and contact
equivalence (change of coordinates and multiplication with a unit or,
equivalently, preserving the isomorphism class of space germs).   It turns out
that right equivalence, which is really a classification of functions, is
easier to handle.   We prove the existence of a finite group $E_{f_0}$
acting on the affine space $T_-$, the base space of the semiuniversal
$\mu$--constant
deformation of $f_0$ of strictly negative weight, such that $T_-/E_{f_0}$ is
the
desired coarse moduli space.   We also show that a fine moduli space almost
never exists.  See \S 1 for definitions and precise statements.
Hence, $T_-/E_{f_0}$ classifies, up to right equivalence, semiquasihomogeneous
power series with fixed principal part.

An important step in the construction of moduli spaces with respect to right
equivalence as well as with respect to contact equivalence is to prove that
isomorphisms between two semiquasihomogeneous
functions have necessarily non--negative degree.   This is proved in \S 2 and
uses the fact that the filtration on the Brieskorn lattice $H''_0(f)$ induced
by the weights coincides with the $V$--filtration, which is independent of the
coordinates.   The proof relies on an analysis of this filtration given in
\cite{He}.

In order to obtain a moduli space with respect to contact equivalence we have
to fix, in addition to the Milnor number, also the Tjurina number.   This is
clear because the dimensions of the
orbits of the contact group acting on $T_-$ depend on the
Tjurina number.   But fixing the Tjurina number is not sufficient.   The orbit
space of the contact group for fixed Tjurina number is, as a topological
space, in general not separated, hence, cannot carry the structure of a
complex space.   It turns out, however, that if we fix the whole Hilbert
function of the Tjurina algebra induced by the weights, the orbit space is a
complex space and a coarse moduli space which classifies, up to contact
equivalence, semiquasihomogeneous hypersurface singularities with fixed
principal part and fixed Hilbert function of the Tjurina algebra.   For
precise statements see \S 4.   These moduli spaces are actually locally closed
algebraic varieties in a weighted projective space.

The orbits of the contact group acting on $T_-$ can also be described as
orbits of an algebraic group $G = U \rtimes (E_{f_0} \cdot \C^\ast)$ where
$E_{f_0}$
is the finite group mentioned above and $U$ is a unipotent algebraic group.
The main ingredient for the proof in the case of contact equivalence is the
theorem on the existence of geometric quotients for unipotent groups in
\cite{GP 2}.   But, in order to give the above simple description of the
strata, we have to use, in a non--trivial way, also the symmetry of the Milnor
algebra, a fact which was already noticed in \cite{LP}.

The stratification with respect to the Hilbert function of the Tjurina algebra
and the proof for the existence of a geometric quotient are constructive and
allow the explicit determination of the moduli spaces and families of normal
forms for specific examples.
\newpage

\section{Moduli spaces with respect to right equivalence}

Let $\C\{x_1, \ldots, x_n\} = \C \{x\}$ be the convergent power series
ring.   Two power series $f, g \in \C\{x\}$ are called {\bf right equivalent}
$(\buildrel r\over\sim)$ if there exists a $\psi \in$ Aut$(\C\{x\})$ such
that $f = \psi(g)$; $f$ and $g$ are called {\bf contact equivalent}
$(\buildrel c\over\sim)$ if there exists a $\psi \in$ Aut$(\C\{x\})$ and $u \in
\C\{x\}^\ast$ such that $f = u \psi(g)$. (Equivalently, the local algebras
$\C\{x\}/(f)$ and $\C\{x\}/(g)$ are isomorphic
respectively the complex germs $(X,0) \subset (\C^n,0)$ and $(Y,0) \subset
(\C^n,0)$ defined by $f$ and $g$ are isomorphic.)

Let $d$ and $w_1, \ldots, w_n$ be any integers.   A polynomial $f_0 \in
\C[x_1, \ldots, x_n] = \C[x]$ is {\bf quasihomogeneous} of {\bf type} $(d;
w_1, \ldots, w_n)$ if for any monomial $x^\alpha = x_1^{\alpha_1} \cdot \ldots
\cdot  x_n^{\alpha_n}$ occurring in $f_0$,

\[
\deg\, x^\alpha := |\alpha| := w_1 \alpha_1 + \cdots + w_n\alpha_n
\]
is equal to $d$.   $w_1, \ldots, w_n$ are called {\bf weights} and $\deg\,
x^\alpha$
is called the (weighted) {\bf degree} of $x^\alpha$.

For an arbitrary power series $f = \sum c_\alpha x^\alpha,\; f \not= 0$, we set
\[
\deg\, f = \inf \{|\alpha|\, \mid\, c_\alpha \not= 0\},
\]
and call it the degree of $f$. For a family of power series $F = \sum
c_{\alpha,
\beta} x^\alpha
s^\beta \in \C \{x,s\}$, parametrized by $\C\{s\}$, we put $\deg_x F= \inf
\{|\alpha| \mid \exists\, \beta$
such that $c_{\alpha,\beta} \not= 0\}$.

$f$ is called quasihomogeneous if it is a
quasihomogeneous polynomial (of some type).   $f$ is called {\bf
semiquasihomogeneous} of type $(d; w_1, \ldots, w_n)$, if
\[
f = f_0 + f_1,
\]
where $f_0$ is a  quasihomogeneous polynomial of type $(d; w_1, \ldots, w_n)$,
$f_1$ is a power series such that $\deg\, f_1 > \deg\, f_0$ and, moreover,
$f_0$ has an isolated singularity at the origin.   $f_0$ is called the {\bf
principal part} of $f$.
Two right equivalent semiquasihomogeneous power series of the same type have
right equivalent principal parts.

Recall (\cite{SaK 1}) that a power series $f$ with isolated singularity is
right equivalent to a quasihomogeneous polynomial
with respect to positive weights if and only if
\[
f \in j(f) := (\partial f/\partial x_1, \ldots, \partial f / \partial x_n).
\]
Moreover, in this case the {\bf normalized weights} $\overline{w}_i =
\frac{w_i}{d}
\in \Q \cap (0,\frac{1}{2}]$ are uniquely determined.

We may consider $f \in \C\{x\}, f(0) = 0$ as a map germ $f : (\C^n,0) \to
(\C,0)$.
An {\bf unfolding} of $f$ over a complex germ or a pointed complex space
$(S,0)$ is by definition a cartesian diagram
\[
\begin{array}{ccc}
(\C^n,0) & \hookrightarrow & (\C^n,0) \times (S,0)\\[1.0ex]
f\; \downarrow & & \downarrow \; \phi\\[1.0ex]
(\C,0) & \hookrightarrow & (\C,0) \times (S,0)\\[1.0ex]
\downarrow & & \downarrow \\[1.0ex]
 0& \hookrightarrow & (S,0).
\end{array}
\]
Hence, $\phi(x,s) = (F(x,s),s)$ and the unfolding
$\phi$ is determined by $F : (\C^n,0) \times (S,0) \to (\C,0),\; F(x,s) = f(x)
+
g(x,s),\; g(x,0) = 0$, and we say that $F$ defines an unfolding of $f$.   Two
unfoldings
$\phi$ and $\phi'$ defined by $F$ and $F'$ over $(S,0)$ are called right
equivalent if
there is an isomorphism $\Psi :(\C^n,0) \times (S,0) \buildrel \cong\over\to
(\C^n,0) \times (S,0),\; \Psi(x,s) = (\psi(x,s),s)$, such that $\phi \circ \Psi
= \phi'$.

For the construction of moduli spaces we have to consider, more generally,
families of unfoldings over arbitrary complex base spaces.   Let $S$ denote a
{\bf category of base spaces}, for example the category of complex germs or of
pointed complex spaces or of complex spaces.   A {\bf family of
unfoldings} over $S \in \ks$ is a commutative diagram
\[
\begin{array}{rcl}
(\C^n,0) \times S & \buildrel\phi\over\longrightarrow & (\C,0) \times S\\
\searrow &  & \swarrow\\
& S &.
\end{array}
\]

Hence, $\phi (x,s) = (G(x,s),s) = (G_s(x),s)$ and for each $s \in S$, the germ
$\phi : (\C^n,0) \times (S,s) \to (\C,0) \times (S,s)$ is an unfolding of $G_s
:
(\C^n,0) \to (\C,0)$.   A morphism of two families of unfoldings $\phi$ and
$\phi' = (G', id_s)$ over $S$ is a morphism $\Psi : (\C^n,0) \times S \to
(\C^n,0) \times S,\; \Psi (x,s) = (\psi(x,s),s) = (\psi_s(x), s)$ such that
$\phi \circ \Psi = \phi'$ (equivalently : $G_s(\psi(x,s)) = G'_s(x)$).
$\phi$ and $\phi'$ are called {\bf right equivalent families of unfoldings} if
there is a morphism $\Psi$ of $\phi$ and $\phi'$ such that for each fixed $s
\in S$, $\psi_s \in \mbox{ Aut}(\C^n,0)$.\\

{}From now on let $f_0  \in \C[x_1, \ldots, x_n]$ denote a quasihomogeneous
polynomial with isolated singularity of type $(d; w_1, \ldots, w_n)$ with $w_i
> 0$ for $i = 1, \ldots, n$.

Consider a power series $f$ which is right equivalent to a
semiquasihomogeneous power series $f'$ of type $(d; w_1, \ldots, w_n)$.
We say that an unfolding $F$ defines an {\bf unfolding of
$\bf f$ of negative weight} over $(S,0)$ if $F$ is right equivalent to $f'(x)
+ g(x,s)$ with $g(x,0) = 0$ and  $\deg_x g > d$.   This holds, for
instance, if there exists a $\C^\ast$--action
with (strictly) negative weights on $(S,0)$ such that $\deg\, g = d$, with
respect to the $\C^\ast$--actions on $(\C^n,0)$ and on $(S,0)$.   By Theorem
2.1 the definition is independent of the choice of $f'$.

We shall now describe the semiuniversal unfolding of $f_0$ of negative weight.
Let $x^\alpha$, $\alpha \in B \subset \N^n$, be a monomial basis of the Milnor
algebra $\C\{x\}/(\partial f_0/\partial x_1, \ldots, \partial f_0/\partial
x_n)$
which is of
$\C$--dimension $\mu$ (the Milnor number of $f_0$), and let $\bar{F} (x,t) =
f_0(x) + \sum_{\alpha \in B} x^\alpha s_\alpha,\; s = (s_\alpha)_{\alpha
\in B} \in \C^\mu$ be the semiuniversal unfolding of $f_0$.   We are mainly
interested in the sub--unfolding  over the affine pointed space $T_- =
(\C^k,0)$,
\[
F(x,t) = f_0(x) + \sum^k_{i=1} t_i m_i,\;\;\; t = (t_1, \ldots, t_k) \in T_-,
\]
where the $m_i$ are the ``upper'' monomials, that is
\[
\{m_1, \ldots, m_k\} = \{ x^\alpha \,\mid\, \alpha\in B,\, |\alpha| > d\}.
\]

For fixed $t \in T_-, F_t(x) = F(x,t) \in \C[x]$ is a semiquasihomogeneous
polynomial with principal part $f_0$.

Let $A = \C[(s_\alpha)_{\alpha \in B}]$ and $A_- = \C[t_1, \ldots, t_k]$.   If
we give weights to $s_\alpha$ and $t_i$ by $w(s_\alpha) = d- |\alpha|$ and
$w(t_i) = d - \deg(m_i)$, then $\bar{F}$ and $F$ are quasihomogeneous
polynomials in $\C[x,s]$ respectively $\C[x,t]$ and $F$ is the restriction of
$\bar{F}$ to $T_-$, the negative weight part of $T =$ Spec\,$A$, defined by
$\{t_1, \ldots, t_k\} = \{s_\alpha\mid
w(s_\alpha) < 0\}$.

{\bf Example}:
$f_0 = x^3 + y^3 + z^7$ is quasihomogeneous of type $(d; w_1, w_2, w_3) = (21;
7, 7, 3)$ with Milnor number $\mu = 24$.   The upper monomials of a monomial
basis of the Milnor algebra $\C\{x,y,z\}/(x^2, y^2, z^6)$ are $m_1 = x z^5,\;
m_2 = yz^5,\; m_3 = xyz^3,\; m_4= xyz^4,\; m_5 = xyz^5$ and,  hence, $A_- =
\C[t_1, \ldots, t_5],\; T_- = \C^5$,
\[
F(x,y,z,t) = f_0 +
\sum^5_{i=1}t_i m_i = f_0 + t_1xz^5 + t_2yz^5 + t_3xyz^3  + t_4xyz^4 +
t_5xyz^5,
\]
$w(t_1, \ldots, t_5) = (-1, -1, -2, -5, -8)$.

\begin{remark}\label{1.1}{\rm
Fix any $t \in T_-$.   $F$ defines an unfolding of $F_t$ of negative weight
over the pointed space $(T_-, t)$.   If we restrict this unfolding to the
germ $(T_-,t)$ this is actually a semiuniversal unfolding of $F_t$ of
negative weight because of the following:

The monomials $m_1, \ldots, m_k$ represent certainly a basis of
$\C \{x\} /(\frac{\partial F_t}{\partial x_1}, \ldots, \frac{\partial F_t}
{\partial x_n})$
for $t$ sufficiently close to $0$, since $\mu(F_t) = \mu(f_0)$.   But, using
the $\C^\ast$--actions on $T_-$ and on $\C^n$, we see that any $F_t$ is contact
equivalent to
some $F_{t'}$, $t'$ close to $0$.
Hence, $\ko_{\C^n \times T_-,0\times T_-}/(\frac{\partial
F}{\partial x_1}, \ldots, \frac{\partial F}{\partial x_n})$ is actually free
over
$T_-$ with basis $m_1, \ldots, m_k$ and the result follows.

We call the affine family
\[
F : \C^n \times T_- \to \C,
\]
$(x,t) \mapsto f_0(x) + \sum^k_{i=1} t_i m_i$ the {\bf semiuniversal family of
unfoldings of negative weight of semiquasihomogeneous
power series with fixed principal part} $\bf f_0$.
}
\end{remark}

\begin{lemma}\label{1.2}
The family of unfoldings $F$ has the following property.  If $f$ is any
semiquasihomogeneous power series with principal part $f_0$, then:

\begin{enumerate}
\item[(i)] $T_- = \{0\}$ if and only if $f_0$ is simple or simple elliptic.
\item[(ii)] There exists a $t \in T_-$ such that $f \buildrel r\over\sim F_t$.
\item[(iii)] Let $f \buildrel r\over\sim F_t$ and let $G(x,s) = f(x) + g(x,s)$
be any
unfolding of $f$ of negative weight over the germ $(S,0)$.   Then there exists
a  morphism, unique on the tangent level, of germs $\varphi : (S,0) \to (T_-,
t)$ such that
$\varphi^\ast F$ is right equivalent to $G$ (that is $T_-$ does not contain
trivial subfamilies of unfoldings).
\item[(iv)] Assume $f_0$ is neither simple nor simple elliptic.     There
exist $t, t' \in T_- ,  t\neq t'$, arbitrarily close to $0$,
such that $F_t \buildrel
r\over\sim F_t'$ (that is $F$ is not universal in any neighbourhood of $0 \in
T_-$).
\end{enumerate}
\end{lemma}
\newpage

{\bf Proof}: \begin{enumerate}
\item[(i)] is due to Saito \cite{SaK 2}.
\item[(ii)] follows from \cite{AGV}, 12.6, Theorem (p.\ 209).
\item[(iii)]  If $T_-$ would contain trivial subfamilies of unfoldings there
must be
a $t \in T_-$ with $\mu(F_t) < \mu(f_0)$, which is not the case.
\item[(iv)] The group $\mu_d$ of $d$--th roots of unit acts on
$T_-$, has $0$ as fixed point and a non--trivial orbit for any $t \not= 0$.
Since for $\xi \in \mu_d,\; F_{\xi\circ t}(\xi \circ x) = \xi^d F_t(x) =
F_t(x)$, two different points of an orbit of $\mu_\alpha$ correspond to right
equivalent functions, we obtain (iv).
\end{enumerate}

Let us introduce the notion of a fine and coarse moduli space for unfoldings of
negative weight with principal part $f_0$ (the weights $w_1, \ldots, w_n$ and
$f_0$ are given as above):  let $\ks$ be a category of base spaces.
   For $S \in \ks$, a {\bf family of unfoldings of negative weight
with principal part} $\bf f_0$ over $S$ is a family of unfoldings
\[
\phi : (\C^n,0) \times S \to (\C,0)\times S,\; (x,s) \mapsto (G(x,s),s)
= (G_s(x),s)
\]
such that: for any $s \in S$,  $G_s :(\C^n,0) \to (\C,0)$ is right equivalent
to a semiquasihomogeneous power series with principal part $f_0$ and the germ
of $G$ at $s$, $G :(\C^n,0) \times (S,s) \to (\C,0)$, is an unfolding of $G_s$
of negative weight. For any
morphism of base spaces $\varphi : T \to S$, the induced map $\varphi^\ast
\phi : (\C^n,0) \times T \to (\C,0)\times T,\; (x,t) \mapsto
(G(x,\varphi(t)),t)$, is an
unfolding of negative weight with principal part $f_0$ over $T$.   Hence, we
obtain a
functor
\[
\mbox{Unf}^-_{f_0} : \ks \to \mbox{ sets}
\]
which associates to $S \in \ks$ the set of right equivalence classes of
families of
unfoldings
of negative weight with principal part $f_0$ over $S$.
If $pt \in \ks$ denotes the base space consisting of one reduced point, then

\begin{tabular}{lp{12cm}}
Unf$^-_{f_0}(pt) =$ & $\{$ right equivalence classes of power series
$f \in \C\{x_1, \ldots, x_n\}$  which are right equivalent to a
semiquasihomogeneous power series with principal part  $f_0\}$.
\end{tabular}

A {\bf fine moduli space} for the functor Unf$^-_{f_0}$ consists of a base
space $T$ and a natural transformation of functors
\[
\psi : \mbox{ Unf}^-_{f_0} \to \mbox{ Hom}(-,T)
\]
such that the pair $(T, \psi)$ represents the functor Unf$^-_{f_0}$.

The pair $(T,\psi)$ is a {\bf coarse moduli space} for Unf$^-_{f_0}$ if

\begin{enumerate}
\item[(i)] if $\psi(pt)$ is bijective, and

\item[(ii)] given the solid arrows (natural transformations) in the following
diagram
\vspace{-0.5cm}

\unitlength1cm
\begin{picture}(8,3)
\put(4,0){Hom$(-,T)$}
\put(6.5,0){\line(1,0){0.15}}
\put(6.7,0){\line(1,0){0.15}}
\put(6.9,0){\line(1,0){0.15}}
\put(7.1,0){\line(1,0){0.15}}
\put(7.3,0){\line(1,0){0.15}}
\put(7.5,0){\line(1,0){0.15}}
\put(7.7,-0.1){$>$}
\put(8.5,0){Hom$(-,T')$,}
\put(6.8,1.5){\vector(-3,-2){1.3}}
\put(7.8,1.5){\vector(3,-2){1.3}}
\put(6.8,2){Unf$^-_{f_0}$}
\end{picture}

there exists a unique dotted arrow (natural transformation) making the diagram
commutative. \end{enumerate}

A fine moduli space is, of course, coarse.

The definitions of fine and coarse moduli spaces still depend on the category
of base spaces $\ks$.   If
$\ks$ is the category of complex germs and if $(S,0) \in \ks$, then Hom$((S,0),
T)$
denotes the set of morphisms of germs $(S,0) \to (T,t)$ where $t$ may be any
point of $T$.   In this case, if $(T, \psi)$ is a fine moduli space, given any
$t \in T$, there exists a unique (up to right equivalence)
universal unfolding of negative weight with principal part $f_0$ over the germ
$(T,t)$ which corresponds to id $\in$ Hom$((T,t), (T,t))$.   But we may not
have a universal family over all of $T$.   If $\ks$ is the category of all
complex spaces, the existence of a fine moduli space implies the existence of
a global universal family over $T$.   But we shall see that even for complex
germs as base spaces a fine moduli space may not exist.   A coarse moduli
space, however, does exist even if $\ks$ is the category of all complex spaces.
The reason is that for a coarse moduli space we do not require any kind of
a universal family.

\begin{theorem}\label{1.3}
Let $E_{f_0}$ be the finite group defined in Definition \ref{2.5}, acting on
$T_-$.   The geometric quotient $T_-/E_{f_0}$ is a coarse moduli space for the
functor Unf$^-_{f_0} : \mbox{ complex spaces } \to$ sets.
\end{theorem}

{\bf Proof}:  Since $E_{f_0}$ is finite, and the action is holomorphic, the
geometric quotient $T_-/E_{f_0}$ exists as a complex space.   According to
Theorem 2.1, Proposition \ref{2.3} and Corollary \ref{2.5}, for any
semiquasihomogeneous
power series $f$ with principal part $f_0$ there exists a unique point
$\underline{t} \in T_-/E_{f_0}$ such that if $f_t \buildrel r\over\sim f$,  $t
\in
T_-$ maps to $\underline{t}$.   In this way we obtain a bijection $\psi(pt)$
from the set of right equivalence classes of semiquasihomogeneous power series
with principal part $f_0$ to $T_-$.

Now let $G : (\C^n,0) \times S \to (\C,0)$ define an element of
Unf$^-_{f_0}(S)$ for some complex space $S$.   We may cover $S$ by open sets
$U_i$ such that there exist morphisms $\varphi_i : U_i \to T_-$ with
$\varphi^\ast F \buildrel r\over\sim G|_{U_i}$.   Even if the $\varphi_i$
are not unique, by the properties of a quotient the compositions $U_i
\buildrel \varphi_i\over \to T_- \to T_-/E_{f_0}$ glue together to give a
morphism
$S \to T_-/E_{f_0}$.   This construction is functorial and provides the desired
natural transformation Unf$^-_{f_0} \to \mbox{ Hom}(-,T_-/E_{f_0})$.   This
finishes
the proof of Theorem \ref{1.3} (for further details for construction of moduli
spaces via geometric quotients cf.\ \cite{Ne}).

\begin{remark}{\rm
(i)\quad If $f_0$ is simple or simple elliptic, then the coarse moduli space
constructed above consists of one reduced point.   Hence, it is even a fine
moduli space.

(ii)\quad If $f_0$ is neither simple nor simple elliptic, Unf$^-_{f_0}$ does
not admit a
fine moduli space, even not if we take complex germs as base spaces.   This
can be seen as follows:  assume there exists such a fine moduli space then,
since it is also coarse, it must be isomorphic to $T_-/E_{f_0}$.   Moreover,
there
exists a universal unfolding over the germ $(T_-/E_{f_0}, 0)$ which can be
induced
from the semiuniversal unfolding $F$ over the germ $(T_-,0)$ and vice versa.
Since $T_-$ does not contain trivial subfamilies, the semiuniversal family $F$
over $(T_-,0)$ would be universal, which contradicts Lemma \ref{1.2} (iv).
}
\end{remark}

{\bf Example}:  Let $f_0(x,y) = x^4 + y^5$.   We obtain $T_- = \C$ and
$F(x,y,t) = x^4 + y^5 + tx^2y^3,\; (d;w_1,w_2;w(t)) = (20;4,5;-2)$.
In this case $E_{f_0} = \mu_d$ and the ring of invariant functions on $T_-$ is
$\C[t^{10}]$, hence $T_-/E_{f_0} \cong \C$.   We give a computational argument
that a fine moduli space does not exist:\\
A local universal family over $(T_-/E_{f_0},0)$ would be given by $G :(\C^n,0)
\times (T_-/E_{f_0}) \to (\C,0),\; (x,y,s) \mapsto G(x,y,s)$.   The proof of
Theorem \ref{1.3} shows that then $F$ would be induced from $G$ by the
canonical map $T_- \to T_-/E_{f_0}$, which is not an isomorphism.
Moreover, the fibre $F^{-1}(0)$ would be isomorphic to
$G^{-1}(0)$ under the map $(x,y,t) \mapsto (x,y,s = t^{10})$.   The image of
this map can be computed by eliminating $t$ from $F(x,y,t) = 0,\; s - t^{10} =
0$.   The result is the hypersurface defined by
$G = (x^4 + y^5)^{10} - sx^{20}y^{30}$.
The special fibre for $s = 0$ has a non--isolated singularity, hence is not
isomorphic to $f_0 = 0$.

\begin{remark}\label{1.5}{\rm
Since the group $E_{f_0}$ acts even algebraically on $T_-$ by Proposition
\ref{2.4}, $T_-/E_{f_0}$ is an algebraic variety.   We may take the category of
base spaces $\ks$ to be the category of (separated) algebraic spaces and
define (families of) unfoldings in the same manner as above, replacing the
analytic local ring $\C\{x\}$ by the henselization of $\C[x]$.   With the same
proof as above we obtain that $T_-/E_{f_0}$ is a coarse moduli space for the
functor
\[
\mbox{Unf}^-_{f_0} : \mbox{ algebraic spaces } \to \mbox{ sets.}
\]
}
\end{remark}
\newpage

\section{Isomorphism of semiquasihomogeneous singularities}

We fix weights $w_1,...,w_n \in \N$ and a degree $d\in \N$
such that the normalized weights $\overline{w}_i = {w_i \over d}$
fulfill $0 < \overline{w}_i \leq {1\over 2}$.
The weights induce a filtration on $\C \{x\}$.
An automorphism $\varphi \neq id$ of $\C \{x\} $ has degree
$m=\deg \varphi$ if $m$ is the maximal number such that
\[
\deg (\varphi (x_i) - x_i) \geq w_i + m \ \ \forall \  i=1,...,n.
\]

The automorphisms of degree $\geq 0$ form the group
$\mbox{Aut}_{\geq 0}(\C \{x\} )$ of all automorphisms of
$\C \{x\} $ which respect the filtration.
The automorphisms of degree $>0$ form a normal subgroup
$\mbox{Aut}_{>0}(\C \{x\} )$ in Aut$_{\geq 0}(\C \{x\} )$.
Automorphisms will be called quasihomogeneous if they map each
quasihomogeneous polynomial to a quasihomogeneous polynomial of the same
degree.
They form a group $G_w \subset \mbox{ Aut}_{\geq 0}(\C \{x\} ),$ which is
isomorphic
to the quotient Aut$_{\geq 0}(\C \{x\} )/\mbox{Aut}_{>0}(\C \{x\} )$.

The image $\varphi (f)$ of a semiquasihomogeneous power series $f$ of degree
$d$ by an automorphism $\varphi$ of $\C \{x\} $ is semiquasihomogeneous
of the same degree if $\deg \varphi \geq 0 $.
The converse is true, too:

\begin{theorem} \label{2.1}
Let $f$ and $g$ be semiquasihomogeneous of degree $d$, and let $\varphi$ be an
automorphism of $\C \{x\} $ such that $\varphi (f) = g$.   Then
$\deg \varphi \geq 0$.
\end{theorem}

{\bf Proof:}
The proof uses some facts which come from the Gauss--Manin connection
for isolated hypersurface singularities
(\cite{SS}, \cite{SaM}, \cite{AGVII}, \cite{He}). The main idea is the
following:
in the case of a semiquasihomogeneous singularity the weights $\overline{w}_i$
induce a filtration on $\C \{x\}$
and a filtration on the Brieskorn lattice $H_0''(f)$.
This last filtration coincides with the $V$--filtration and is independent
of the coordinates.

The Brieskorn lattice $H_0''(f)$ is
\[
H_0'' = \Omega^n / df\land d\Omega^{n-1} .
\]
Here $\Omega^k = \Omega^k_{\C^n,0}$ denotes the space of germs of
holomorphic $k$--forms. The class of $\omega \in \Omega^n$ in $H_0''(f)$
is denoted by $s[\omega ]_0 \in H_0''(f)$.
The $V$--filtration on $H_0''(f)$ is determined by the orders
$\alpha_f (\omega) = \alpha_f(s[\omega]_0)$ of $n$--forms
$\omega \in \Omega^n$.
The most explicit description of the order $\alpha_f (\omega)$ might be
the following (\cite{AGVII}, \cite{He}):
\newpage

\begin{eqnarray*}
\alpha_f(\omega) = \min\,\{ \alpha & | & \exists
	\mbox{ (manyvalued) continuous family of cycles}\\
  &&\delta (t)\in H_{n-1}(X_t,\Z) \mbox{ on the Milnor fibers }X_t \\
  &&\mbox{of the singularity }f:(\C^n,0) \to (\C,0),\\
  &&\mbox{such that } a_{\alpha,k} \neq 0 \mbox{ in } \\
  &&\int\limits_{\delta (t)} {\omega \over df} =
       \sum_{\beta,k} a_{\beta,k}\cdot t^\beta \cdot (\ln t)^k \\
  &&\mbox{for a } k \mbox{ with } 0\leq k \leq n-1\ \}.
\end{eqnarray*}

The description shows that we have
\[
\alpha_f (\omega) = \alpha_g (\varphi (\omega))
  = \alpha_g (\varphi(h) d\varphi(x))
\]
for $\omega = h(x)dx_1...dx_n = hdx \in \Omega^n$.

Since $f$ is semiquasihomogeneous it is possible to give a simple algebraic
description of the order $\alpha_f (\omega)$.
Indeed, we define mappings
\vspace{-0.5cm}

\begin{eqnarray*}
  \nu_C & : & \C\{ x_1,...,x_n\} \to \Q_{\geq 0}\cup \{\infty\},\\
  \nu_\Omega & : & \Omega^n \to \Q_{>-1}\cup \{\infty\},\\
  \nu_f & : & H_0''(f) \to \Q_{>-1}\cup \{\infty\}
\end{eqnarray*}

by
\[ \nu_C (x^\alpha) = \sum_{i=1}^n \overline{w}_i \alpha_i\ , \
    \nu_C(0) = \infty,  \
    \nu_C (\sum b_\alpha x^\alpha) =
           \min \{\nu_C (x^\alpha)\ |\ b_\alpha\neq 0\}
\]
and
\[
\nu_\Omega (hdx) = \nu_C (hx_1...x_n) -1
\]
and
\[
\nu_f (s[\omega]_0) = \nu_f(\omega) =
    \max \{\nu_\Omega (\eta)\ |\ s[\eta]_0 = s[\omega]_0\}.
\]

Then, from \cite{He}, Chapter 2.4, it follows that
\[
\nu_f (\omega) = \alpha_f (\omega) = \alpha_g (\varphi(\omega))
		 = \nu_g (\varphi(\omega)).
\]
For all $\eta\in \Omega^{n-2}$ we have
\[
\nu_f (df\land d\eta)
   \geq -1 + \sum_j \overline{w}_j + (1-\max( \overline{w}_i ))
   \geq \sum_j \overline{w}_j - \frac{1}{2}.
\]
For $\omega$ with
\[
\min \{\nu_\Omega (\omega ), \nu_f (\omega ),\nu_g (\varphi (\omega )),
    \nu_\Omega (\varphi (\omega )) \}
       < \sum_j \overline{w}_j - \frac{1}{2}
\]
this implies
\[
\nu_\Omega (\omega ) = \nu_f (\omega ) = \nu_g (\varphi (\omega))
       = \nu_\Omega (\varphi (\omega )).
\]
We obtain
\[
\sum_j \overline{w}_j -1 = \nu_\Omega (dx) = \nu_f (dx)
       = \nu_g (d\varphi(x)) = \nu_\Omega (d\varphi(x)).
\]
For $i$ with $\overline{w}_i < \frac{1}{2}$ we obtain
\vspace{-0.5cm}

\begin{eqnarray*}
  \overline{w}_i + \nu_\Omega (dx) & = & \nu_C(x_i) + \nu_\Omega (dx)
    = \nu_\Omega (x_idx) = \nu_f (x_idx) \\
  & = & \nu_g (\varphi (x_i) d\varphi(x))
    =   \nu_\Omega (\varphi (x_i) d\varphi(x))
    =   \nu_C (\varphi (x_i)) + \nu_\Omega (d\varphi(x))  \\
  & = & \nu_C (\varphi (x_i)) + \nu_\Omega (dx)
\end{eqnarray*}
and $ \nu_C (\varphi (x_i)) = \overline{w}_i $.

For $i$ with $\overline{w}_i = \frac{1}{2}$ the equality $\nu_\Omega (x_idx) =
\sum \overline{w}_i - \frac{1}{2} $
implies
\[
\nu_\Omega (\varphi (x_idx)) \geq
  \sum \overline{w}_i - \frac{1}{2}
\]
and $ \nu_C (\varphi (x_i)) \geq \frac{1}{2} $.
Therefore, $\nu_C (\varphi (x_i)) \geq \nu_C (x_i) = w_i \ \ \forall
i=1,...,n,$ and thus $\deg \varphi \geq 0$. \hfill

\begin{remark}\label{2.2} \rm
In the following, Theorem \ref{2.1} will be used to describe
a finite group $E_{f_0}\subset \mbox{ Aut}(T_{-})$ which operates transitively
on each set of parameters in $T_{-}$ which belong to one right equivalence
class.
Theorem \ref{2.1} also shows that the Hilbert function
of the Tjurina algebra (cf.\ Chapter 4) is an invariant of the
contact equivalence class.
\end{remark}

Now let $f_0\in \C[x_1,...,x_n]$ be quasihomogeneous of degree $d$
with an isolated singularity in 0. Let $m_1,...,m_k$ denote the monomials
of degree $>d$ in a monomial base of the Milnor algebra of $f_0$.
Consider the semiuniversal unfolding of $f_0$ of negative weight,
\[
F = f_0 + \sum_{i=1}^k m_it_i.
\]
For a fixed value of $t$ we write $F_t = f_0 + \sum m_it_i$.
With $\deg t_i = w(t_i) = d - \deg m_i < 0$ we obtain a filtration on
$\C[t_1,...,t_k] = A_{-}$ such that $F \in \C[x,t]$ is
quasihomogeneous of degree $d$ in $x$ and $t$. We write $T_{-}=Spec\ A_{-}$
(cf.\ \S 1).

\begin{proposition} \label{2.3}
For any semiquasihomogeneous power series $f$ with principal part $f_0$ there
exist an automorphism $\varphi \in Aut_{>0}(\C\{x\})$ and a
parameter $t\in T_{-}$ such that $\varphi (f) = F_t$.
The $t\in T_{-}$ is uniquely determined.
\end{proposition}

{\bf Proof}:
The existence of $\varphi$ and $t$ is proved in \cite{AGV}, 12.6,
Theorem (p.\ 209).
The following proves the uniqueness of $t$.

Let $t$ and $t'\in T_{-}$ and $\psi \in \mbox{ Aut}_{>0}(\C\{x\})$ be given
such
that $\psi (F_t) = F_{t'}$.
With $\psi_s(x_i) = x_i + s(\psi(x_i)-x_i)$ we obtain a family $\psi_s$ of
automorphisms in Aut$_{>0}(\C\{x\})$.
The family $\psi_s(F_t)$ of semiquasihomogeneous functions with principal
part $f_0$ connects $\psi_0(F_t) = F_t$ and $\psi_1(F_t) = F_{t'}$.
The family may not be contained in $T_{-},$ but can be induced from $T_{-}$
by a suitable base change:
Following the proof of the theorem in [AGV], 12.6 (p. 209), we can find a
family $\chi_s$ of automorphisms and a holomorphic map $\sigma: \C \to T_-$
such that $\chi_s \circ \psi_s (F_t) = F_{\sigma (s)}$ and
$\chi_s\in Aut_{>0}(\C \{x\})$ and even $\chi_0 = id = \chi_1,\
\sigma (0)=t,\, \sigma (1)=t'$.
But since $T_{-}$ is part of the
semiuniversal deformation,
which is miniversal on the $\mu$--constant stratum,
and since $T_{-}$ does not contain trivial subfamilies
with respect to right equivalence, $t=t'$ as desired. \hfill

\begin{proposition}  \label{2.4}
\begin{enumerate}
\item For any $\varphi \in G_w^{f_0} = \{\psi\in G_w\ |\ \psi(f_0)=f_0\}$
      and any $t\in T_{-}$ there exist $s = \theta(\varphi)(t) \in T_{-}$
      and an automorphism $\psi\in Aut_{>0}(\C\{x\})$
      such that $\psi \circ \varphi (F_t) = F_s$.

\item The function $\theta (\varphi) : T_{-} \to T_{-}$ is uniquely
      determined, bijective and fulfills
      $\theta(\varphi^{-1}) = \theta^{-1}(\varphi)$ and
      $\theta(\varphi)\circ \theta(\psi)
	 =\theta(\varphi \circ \psi)$ for any $\psi \in G_w^{f_0}$.

\item The components $\theta (\varphi)(t_i)$ are quasihomogeneous
      polynomials in $A_{-}$ of degree $\deg (t_i)$.
\end{enumerate}
\end{proposition}

{\bf Proof}:
The statements 1.\ and 2.\ follow from Proposition \ref{2.3}
and from the fact that
Aut$_{>0}(\C\{x\})$ is a normal subgroup of
Aut$_{\geq 0}(\C\{x\})$. Statement 3.\ follows from the proof of the theorem
in \cite{AGV}, 12.6 (p.\ 209).
Along the lines of this proof one can construct power series
$\psi_1,...,\psi_n \in \C\{x,t\}$
    and a family of automorphisms $\psi (t)$ such that
    $\psi (t)(x_i) = \psi_i(t)$
with the following properties:

\begin{quote}
$\psi_i$ is quasihomogeneous in $x$ and $t$ of degree $w_i,$  \\
$\psi_i - x_i$ has degree $>w_i$ in $x,$                      \\
for any fixed $t$ the automorphism $\psi (t)\in
	    \mbox{ Aut}_{>0}(\C\{x\})$
with $\psi (t)(x_i) = \psi_i(t)$ gives
$\psi (t) \circ \varphi (F_t) = F_{\theta (\varphi)(t)}$.
\end{quote}

The power series $F=f_0+\sum m_it_i,$ and
$\varphi (F)=f_0+\ldots$ and
$\psi (t)\circ \varphi (F) = f_0 + \sum m_i\theta (\varphi)(t_i)$
are all quasihomogeneous of degree $d$
with respect to $x$ and $t$.
This proves 3. \hfill \\

The functions $\theta (\varphi)$ are biholomorphic.

\begin{definition}
The image
$\theta (G_w^{f_0})$ in Aut$(T_{-})$ will be denoted by $E_{f_0}$.
\end{definition}

\begin{corollary}\label{2.5}
The map $\theta : G_w^{f_0} \to E_{f_0} \subset Aut (T_{-})$
is a group homomorphism.
The automorphisms $\theta (\varphi)$ of $T_{-}$ commute with the
$\C^{*}$-operation on $T_{-}$.
Each orbit of $E_{f_0}$ consists of all parameters in $T_{-}$ which belong to
one right equivalence class.
\end{corollary}

{\bf Proof}:
The first two statements follow from Proposition \ref{2.4},
the third statement follows from Theorem \ref{2.1}. \hfill

\begin{proposition}\label{2.6}
\begin{enumerate}
\item The group $G_w^{f_0}$ is finite if
      $\overline{w}_1,...,\overline{w}_{n-1} < \frac{1}{2}$ and
      $\overline{w}_n \leq \frac{1}{2}$.

\item The group $E_{f_0}$ is finite.
\end{enumerate}
\end{proposition}

{\bf Proof}:\\
{\bf 1.}\quad The dimension of the algebraic group $G_w$ is
\[
\dim G_w = \sum_{i=1}^n \# (\mbox{ monomials } x^\alpha
    \mbox{ of degree }w_i \ ).
\]
The group $G_w$ operates on
\[
V = \bigoplus_{\deg x^\alpha = d} \C\cdot x^\alpha.
\]

Let $j(f_0)$ denote the Jacobi ideal of $f_0$ and $j_i(f_0)$ the ideal
\[
j_i(f_0) = (\frac{\partial f_0}{\partial x_{1}},\ldots,
\frac{\partial f_0}{\partial x_{i-1}},
\frac{\partial f_0}{\partial x_{i+1}},\ldots,
\frac{\partial f_0}{\partial x_{n}}).
\]
The tangent space $T_{f_0}G_w f_0
\subset T_{f_0}V$ of $G_w f_0$ in $f_0$ is
\[
T_{f_0} G_w f_0 \cong j(f_0)\cap V.
\]

For any relation
\[
0 = \sum_{i=1}^n \sum_{\deg x^\alpha = w_i}
	   a_{\alpha,i}\cdot x^\alpha \cdot
	   \frac{\partial f_0}{\partial x_{i}}
     = \sum_{i=1}^n b_i \frac{\partial f_0}{\partial x_{i}}
\]
with $a_{\alpha,i}\in \C$ and $b_i = \sum_{\deg x^\alpha = w_i} a_{\alpha,i}
\cdot x^\alpha $ we have $\deg b_i = w_i$ and
$\deg \frac{\partial f_0}{\partial x_{j}} = d - w_j > w_i$ for $j\neq i$.
Therefore, $b_i \not\in j_i(f_0)$ or $b_i = 0$.
But since $f_0$ has an isolated singularity, the sequence
$(\frac{\partial f_0}{\partial x_{1}},\ldots, \frac{\partial f_0}{\partial
x_{n}}) $ is a regular sequence and $\frac{\partial f_0}{\partial x_{i}}$
is not a zero divisor in $j_i(f_0)$. This implies $b_i = 0 $ for any $i,$ and
\[
j(f_0)\cap V = \bigoplus_{i=1}^n\;
	 \bigoplus_{\deg x^\alpha = w_i} \C \cdot x^\alpha \cdot
	 {\partial f_0 \over \partial x_{i}},
\]
and
\[
\dim G_w^{f_0} = \dim G_w - \dim j(f_0) \cap V = 0.
\]

{\bf 2.}\quad One can order the weights $w_i$ such that
$\overline{w}_1,\ldots,\overline{w}_r < \frac{1}{2},$
$\overline{w}_{r+1},\ldots,\overline{w}_n = \frac{1}{2}$.
The generalized Morse lemma and Theorem \ref{2.1} imply
the existence of an automorphism $\varphi \in G_w$ and of a
quasihomogeneous polynomial $g_0 \in \C[x_1,\ldots,x_r]$ of degree $d$
such that
$\varphi (f_0) = g_0 + x^2_{r+1} + \ldots + x^2_n $.
Now let $\widetilde{m}_1,\ldots,\widetilde{m}_k$ be the monomials of degree $>
d$
in a
monomial base of the Jacobi algebra of $g_0$.
Analogously to $F$ we obtain families
\vspace{-0.5cm}

\begin{eqnarray*}
     \widetilde{G} & = & g_0 + \sum_{i=1}^k
			       \widetilde{m}_i \widetilde{t}_i \\
     \mbox{and } G & = & \widetilde{G} + x^2_{r+1} + \ldots + x^2_n .
\end{eqnarray*}

It is well known that $G_t$ and $G_{t'}$ are right equivalent
if and only if $\widetilde{G}_t$ and $\widetilde{G}_{t'}$ are right equivalent.
Let
$\widetilde{w}$ be the tuple of weights $\widetilde{w} = (w_1,\ldots,w_r)$.
The group $G_{\widetilde{w}}^{g_0}$ is finite by
the first part of this proposition and induces a finite
group $\widetilde{E}_{\widetilde{w}}$ of automorphisms of
$\widetilde{T}_{-} = Spec\ \C[\tilde{t}]$.
In fact this is the largest subgroup of Aut$(\widetilde{T}_{-})$ which respects
the right equivalence classes.
Similarly to Proposition \ref{2.4} one can prove that $\varphi$
induces a biholomorphic mapping from $T_{-}$ to $\widetilde{T}_{-}$
which respects the right equivalence classes.
This gives an injective (in fact bijective) mapping from $E_{f_0}$ to
$\widetilde{E}_{\widetilde{w}}$.
Hence, $E_{f_0}$ is finite. \hfill

\begin{example}\label{2.7} \rm
$f_0 = x^3 + y^3 + z^7,\ (d;w_1,w_2,w_3) = (21;7,7,3),\ T_{-} = \C^5$,
$F = f_0 + \sum_{i=1}^5 t_im_i = f_0 + t_1xz^5 + t_2yz^5 + + t_3xyz^3 +
t_4xyz^4
+ t_5xyz^5$,
the weights of $(t_1,\ldots,t_5)$ are $(-1,-1,-2,-5,-8)$.

$G_w^{f_0}$ contains $6\cdot 3 \cdot 7$ elements:
obviously, $G_w^{f_0} \cong G_{(1,1)}^{g_0} \times \Z_7$
where $g_0 = x^3 + y^3$.
The group $G_{(1,1)}^{g_0}$ is isomorphic to a subgroup of
${\bf Gl}(2,\C)$.
The image in ${\bf PGl}(2,\C)$ permutes three points in ${\bf P^1C}$
and is isomorphic to $S_3,$ the kernel is isomorphic to
$\{ id , \xi \cdot id , \xi^2 \cdot id\},$ where $\xi = e^{2\pi i/3}$.
Therefore $G_{(1,1)}^{g_0}$ is
\[
G_{(1,1)}^{g_0} = (\langle\alpha\rangle \ltimes \langle\beta\rangle)
			  \times \langle\gamma\rangle
			  \times \langle\delta\rangle
   \cong S_3 \times \Z_3 \times \Z_7
\]
with
\vspace{-0.5cm}

\begin{eqnarray*}
   \alpha & : \ (x,y,z)\  \to & (    y,      x,z),  \\
   \beta  & : \ (x,y,z)\  \to & (\xi x,\xi^2 y,z),  \\
   \gamma & : \ (x,y,z)\  \to & (\xi x,\xi   y,z),  \\
   \delta & : \ (x,y,z)\  \to & (    x,      y,e^{2\pi i/7}z).
\end{eqnarray*}

The mapping $\theta: G_w^{f_0} \to E_{f_0}$ is an isomorphism with
\vspace{-0.5cm}

\begin{eqnarray*}
   \theta(\alpha) & :\  (t_1,t_2,t_3,t_4,t_5)\  \to &
       ( t_2, t_1, t_3, t_4 , t_5) ,\\
   \theta(\beta ) & :\  (t_1,t_2,t_3,t_4,t_5)\  \to &
       (\xi t_1,\xi^2 t_2, t_3, t_4, t_5) ,\\
   \theta(\gamma) & :\  (t_1,t_2,t_3,t_4,t_5)\  \to &
       (\xi t_1,\xi t_2,\xi^2 t_3,\xi^2 t_4 ,\xi^2 t_5) ,\\
   \theta(\delta) & :\  (t_1,t_2,t_3,t_4,t_5)\  \to &
       (\zeta^5 t_1,\zeta^5 t_2,\zeta^3 t_3,\zeta^4 t_4 ,\zeta^5 t_5)
       \mbox{ with } \zeta = e^{2\pi i/7}.
\end{eqnarray*}

Let $\C^{\ast}$ denote the group of $\C^{\ast}$-operations on $T_{-}$.
Then $E_{f_0} \cap \C^{\ast} = \langle \theta(\gamma),\theta(\delta)\rangle$
and  $E_{f_0}\cdot \C^{\ast} \cong \langle\theta(\alpha),\theta(\beta)\rangle
\times \C^{\ast} \cong S_3\times\C^{\ast}$.
\end{example}
\newpage

\section{Kodaira--Spencer map and integral manifolds}

Let $f_0$ be semiquasihomogeneous of type $(d; w_1, \ldots, w_n)$, $w_i > 0$,
and $F : \C^n \times T_- \to \C,\; (x,t) \mapsto f_0(x) + \sum\limits^k_{i=1}
t_im_i$,\enspace the semiuniversal family of unfoldings of negative weight as
in \S
1.   In order to describe the orbits of the contact group acting on $T_-$ we
study the Kodaira--Spencer map of the induced semiuniversal family of
deformations (of space germs) defined as follows.   Let
\[
\kx = \{(x,t) \in \C^n \times T_- \mid F(x,t) = 0\}
\]
and let $(\kx, 0 \times T_-)$ denote the germ of $\kx$ along the trivial
section $0 \times T_-$ which is a subgerm of
$(\C^n \times T_- , 0 \times T_-) = (\C^n,0) \times T_-$.
The composition with the projection gives a morphism
\[
\phi : (\kx, 0 \times T_-) \hookrightarrow (\C^n,0) \times T_- \to T_-
\]
such that, for any $t \in T_-$, $(\phi^{-1}(t), (0,t)) \cong (\kx_t,0)\subset
(\C^n,0)$ is a semiquasihomogeneous hypersurface singularity with principal
part equal to $(\kx_0,0) = (f_0^{-1} (0),0) =: (X_0,0)$.    We call this family
the {\bf semiuniversal family of deformations of negative weight of
semiquasihomogeneous hypersurface singularities with fixed principal part}
$(\bf X_0,0)$ (see also \S 4).

For the study of the Kodaira--Spencer map of $(\kx,0 \times T_-) \to T_-$ it
is more convenient to work on the ring level $A_- \to A_-\{x\}/F$.

The Kodaira--Spencer map (cf.\ \cite{LP}) of the family $A_- \to A_-\{x\}/F$,
\[
\rho : \mbox{Der}_\C A_- \to (x)A_-\{x\}/\left(F + (x) (\frac{\partial
F}{\partial x_1}, \ldots, \frac{\partial F}{\partial x_n})\right),
\]
is defined by $\rho(\delta) = \mbox{ class}(\delta F) = \mbox{
class}(\sum\limits^k_{i=1} \delta(t_i)m_i)$.

Let $\kl$ be the kernel of $\rho$.   $\kl$ is a Lie--algebra and along the
integral manifolds of $\kl$ the family is analytically trivial (cf.\
\cite{LP}).

In our situation it is possible to give generators of $\kl$ as $A_-$--module:

Let $I = A_-\{x\}/(\frac{\partial F}{\partial x_1}, \ldots, \frac{\partial
F}{\partial x_n})$, then $I$ is a free $A_-$--module and $\{m_i\}_{i=1,
\ldots, k}$ can be extended to a free basis.

Multiplication by $F$ defines an endomorphism of $I$ and $F I \subseteq
\bigoplus\limits^k_{i=1} m_iA_-$.

Define $h_{\alpha,j}$ by
\[
x^\alpha F = \sum h_{\alpha,j} m_j \mbox{ in } I.
\]

Then $h_{\alpha j}$ is homogeneous of degree $|\alpha| + \deg(t_j) = |\alpha|+
d -
\deg(m_j)$.   This implies $h_{\alpha j} = 0$ if $|\alpha| + \deg(t_j) \ge 0$,
in
particular $h_{\alpha j} = 0$ if $|\alpha| \ge (n-1)d -2 \sum w_i$.   For
$\alpha$
and $|\alpha| < (n-1)d -2 \sum w_i$ let $\delta_\alpha := \sum
h_{\alpha,j} \frac{\partial}{\partial t_j}$.

\begin{proposition}\label{3.1} (cf.\ \cite{LP}, Proposition 4.5):
\begin{enumerate}
    \item $\delta_{\alpha}$ is homogeneous of degree $|\alpha|$.
    \item $\kl = \sum A_- \delta_\alpha$.
\end{enumerate}
\end{proposition}

Now there is a non--degenerate pairing on $I$ (the residue pairing) which is
defined by $\langle h, k\rangle = \mbox{ hess} (h \cdot k)$.   Here
$hess(h)$ is the evaluation of $h$ at the socle (the hessian of $f$).

Using the pairing one can prove the following:

\begin{proposition}\label{3.2}
There are homogeneous elements $n_1, \ldots, n_k \in A_- \{x\}$ with the
following properties:

\begin{enumerate}
    \item If $n_i F = \sum^k_{j = 1} h_{ij} m_j$ in $I$
    then $h_{ij} = h_{k-j+1, k-i+1}$.
    \item If $\delta_i := \sum^k_{j=1} h_{ij}
    \frac{\partial}{\partial t_j}$ then
    $\delta_i$ is homogeneous of degree $\deg(n_i)$ and
    $\kl =  \sum^k_{i = 1} A_-\delta_i$.
\end{enumerate}
\end{proposition}

In \cite{LP} (Proposition 5.6) this proposition is proved for $n = 2$.   The
proof can easily
be extended to arbitrary $n$.   The important fact is the symmetry, expressed
in 1.\\

Let $L_+$ be the Lie--algebra of all vector fields of $\kl$ of degree $\ge w =
\min\{w_i\}$.   Then
$L$ is finite dimensional and nilpotent.   $\delta_2, \ldots, \delta_k \in
L_+$ and $\delta_1 = \sum\limits^k_{i=1} \deg(t_i) t_i
\frac{\partial}{\partial t_i}$ is the Euler vector field (cf. \cite{LP}).
Let $L = L_+ \oplus \C \delta_1$ then $L$ is a finite dimensional and solvable
Lie--algebra and $\kl = \sum A_- L,\; L/L_+ \cong \C\delta_1$.

\begin{corollary}\label{3.3}
The integral manifolds of $\kl$ coincide with the orbits of the
algebraic group $exp(L)$.
\end{corollary}

Now consider the matrix $M(t) := (\delta_i(t_j))_{i, j = 1,
\ldots, k} = (h_{ij})_{i,j=1, \ldots, k}$. Evaluating this matrix
at $t \in T_-$ we have

\begin{eqnarray*}
\mbox{rank } M(t) & = & \mbox{dimension of a maximal integral manifold of }
\kl\\
		  &   & \mbox{(resp.\ of the orbit of exp(L)) at}\; t\\
		  & = & \mu - \tau(t),
\end{eqnarray*}

where $\tau(t)$ denotes the Tjurina number of the singularity defined by $t$
\enspace i.e.\ of $F(x,t)$.

\begin{example}\label{3.4}{\rm
We continue with Example \ref{2.7}, $f_0 = x^3 + y^3 + z^7$.   Let
\vspace{-0.5cm}

\begin{eqnarray*}
	    n_1 & = & -21\\
	    n_2 & = & -21 z + \left(\frac{250}{49} t_1^3 t_2 + \frac{55}{7}
	    t_1^2 t_3 - \frac{250}{49} t_2^4\right) y - \frac{55}{7} t_2^2 t_3
	    x\\
	    n_3 & = & -21 z^2 - 30 t_2 y\\
	    n_4 & = & -21x\\
	    n_5 & = & -21y
\end{eqnarray*}
then the matrix defined by Proposition \ref{3.2} is
\[
(\delta_i(t_j)) = \left(
\begin{array}{ccccc}
t_1 & t_2 & 2t_3 & 5t_4 & 8t_5\\
0 & 0 & 0 & 2t_3-\frac{10}{7}t_1t_2 & 5t_4\\
0 & 0 & 0 & 0 & 2t_3\\
0 & 0 & 0 & 0 & t_2\\
0 & 0 & 0 & 0 & t_1
\end{array}
\right).
\]

We have $\mu = 24$ and

\begin{tabular}{lp{14cm}}
$\tau =  21$ & if and only if  $2t_3 - \frac{10}{7} t_1t_2 \not= 0$,\\
$\tau =  22$ & if and only if  $2t_3 - \frac{10}{7} t_1t_2 = 0$ and $t_1
\not= 0$ or $t_2 \not= 0$ or $t_3 \not= 0$ or $t_4 \not= 0$,\\
$\tau = 23$ & if and only if $t_1 = t_2 = t_3 = t_4 = 0$ and $t_5 \not= 0$,\\
$\tau = 24$ & if and only if $t_1 = t_2 = t_3 = t_4 = t_5 = 0$.
\end{tabular}}
\end{example}
\newpage

\section{Moduli spaces with respect to contact equivalence}

In this section we want to construct a coarse moduli space for
semiquasihomogeneous hypersurface singularities with fixed principal part with
respect to contact equivalence, that is isomorphism of space germs.
Such a moduli space does only exist if we
fix further numerical invariants.   We shall use the Hilbert function of the
Tjurina algebra induced by the given weights.

Let us first define the functor for which we are going to construct the moduli
space.


A complex germ $(X,0) \subset (\C^n,0)$ is called a {\bf quasihomogeneous}
(respectively {\bf semiquasihomogeneous}) {\bf hypersurface singularity} of
type $(d; w_1, \ldots, w_n)$ if there exists a quasihomogeneous polynomial $f
\in \C[x_1, \ldots, x_n]$ (respectively a semiquasihomogeneous power series $f
\in \C\{x_1, \ldots, x_n\})$ of type $(d; w_1, \ldots, w_n)$ such that $(X,0)
= (f^{-1} (0),0)$.   If $f_0$ is the principal part of $f$ then $(X_0,0) =
(f_0^{-1}(0),0)$ is called the {\bf principal part} of $(X,0)$.   Multiplying
$f$ with a unit changes $f_0$ by a constant, hence the principal part if
well--defined.   Two power series are contact equivalent if and only if the
corresponding space germs are isomorphic.

A {\bf deformation} ({\bf with section}) of $(X,0)$ over a complex germ or a
pointed complex space $(S,0)$ is a cartesian diagram
\[
\begin{array}{ccc}
0 & \hookrightarrow & (S,0)\\
\downarrow & & \downarrow\; \sigma\\
(X,0) & \hookrightarrow & (\kx,0)\\
\downarrow & & \downarrow\; \phi\\
0 & \hookrightarrow & (S,0)
\end{array}
\]
such that $\phi$ is flat and $\phi \circ \sigma =$ id.   Two deformations
$(\phi, \sigma)$ and $(\phi', \sigma')$ of $(X,0)$ over $(S,0)$ are isomorphic
if there is an isomorphism $(\kx,0) \buildrel\cong\over\to (\kx',0)$ such that
the obvious diagram commutes.   We shall only consider deformations with
section.

If $(X,0) = (f^{-1}(0),0)$ and if $F : (\C^n,0) \times (S,0) \to (\C,0)$ is an
unfolding of $f$ then the projection $(\kx,0) = (F^{-1}(0),0) \to (S,0)$ is a
deformation of $(X,0) \hookrightarrow (\kx,0)$ with trivial section $\sigma(s)
= (0,s)$.   Conversely, any deformation of $(X,0)$ is isomorphic to a
deformation induced by an unfolding in this way.   A  deformation
($\phi,\sigma$) of a hypersurface singularity $(X,0)$, which is isomorphic to a
semiquasihomogeneous hypersurface singularity $(X',0) = (f^{-1}(0),0)$ of type
$(d; w_1, \ldots, w_n)$ over $(S,0)$, is called {\bf deformation of negative
weight} if it is isomorphic to a deformation induced by an unfolding of $f$
of negative weight.

We have to show that the definition is independent of the chosen unfolding:
two inducing unfoldings differ by a right equivalence and a multiplication
with a unit.   We have shown in \S 1 that the definition depends only on the
right equivalence class.   Hence, we have to show the following:  if $f(x)$ is
a semiquasihomogeneous power series, $f(x) + g(x,s),\; g(x,0) = 0,\; \deg_x g >
d$, an unfolding of negative weight and $u(x,s) \in \ko^\ast_{\C^n \times
S,0}$ a unit, then $u(f+g) \buildrel r\over\sim f'(x) + g'(x,s)$ with
$f^{-1}(0) = f'^{-1}(0),\; g'(x,0) = 0$ and $\deg_x g'>d$.   Replacing
$u(x,s)$ by $(u(x,0))^{-1} u(x,s)$ we may assume that $u(x,s) = u_0(s) +
su_1(x,s),\; u_0(0) = 1,\; u_1(0,s) = 0$.
If $\nu \in \ko_{S,0}$ is a $d$--th root of $u_0$ and if $\psi$ denotes the
automorphism of degree 0, $\psi (x,s) = (\nu(s)^{w_1}x_1, \ldots, \nu(s)^
{w_n}x_n)$, then $u_0(s) f(x)
= f(\psi(x,s)) + s \tilde{f}(x,s),\; \deg_x \tilde{f} > d$.   But this implies
$u(f+g) \circ \psi^{-1} = f + g'$ with $g'(x,0) = 0$ and $\deg\, g'_x > d$ as
desired.

Again, we have to consider not only germs but also arbitrary complex spaces as
base spaces.   A {\bf family of deformations} of hypersurface singularities
over a base space $S\in \ks$ is a morphism $\phi : \kx \to S$ of complex
spaces together with a section $\sigma : S \to \kx$ such that for each $s \in
S$ the morphism of germs $\phi : (\kx, \sigma(s)) \to (S,s)$ is flat and the
fibre $(\kx_s, \sigma(s)) = (\phi^{-1} (s), \sigma(s))$ is a hypersurface
singularity.   This is, of course, only a condition on the germ $(\kx,
\sigma(S))$ of $\kx$ along $\sigma(S)$.   A morphism of two families
$(\phi,\sigma)$ and $(\phi',\sigma')$ over $S$ is a morphism $\psi : \kx \to
\kx'$ such that $\phi = \phi' \circ \psi$ and $\sigma' = \psi \circ \sigma$.
$(\phi,\sigma)$ and $(\phi',\sigma')$ are called {\bf contact equivalent} or
{\bf isomorphic families of deformations} if there exists a morphism  $\psi$
such that for any $s \in S$, $\psi$ induces an isomorphism of the germs of the
fibres $(\kx_s,\sigma(s)) \cong (\kx'_s, \sigma'(s))$.

Let us fix a quasihomogeneous hypersurface singularity $(X_0,0) \subset
(\C^n,0)$ of type $(d; w_1, \ldots, w_n)$.   For $S \in \ks$, a {\bf family of
deformations of negative weight with
principal part} $(\bf X_0, 0)$ over $S$ is a family of deformations
\[
S \buildrel \sigma\over\to (\kx, \sigma(S)) \buildrel\phi\over\to S
\]
with section such that:  for any $s \in S$ the fibre $(\kx_s, \sigma(s))$ is
isomorphic to a semiquasihomogeneous hypersurface singularity with principal
part $(X_0,0)$ and the germ $(S,s) \buildrel\sigma\over\to (\kx,\sigma(s))
\buildrel \phi\over\to (S,s)$ is a deformation of $(\kx_s,\sigma(s))$ of
negative weight.

For any morphism of base spaces $\varphi : T \to S$, the induced deformation
$T \to (\varphi^\ast \kx,\; \varphi^\ast \sigma(T)) \to T$ is a family of
deformations with negative weight and principal part $(X_0,0)$.   We obtain a
functor
\[
\mbox{Def}^-_{X_0} : \ks \to \mbox{ sets}
\]
which associates to $S \in \ks$ the set of isomorphism classes of families of
deformations of negative weight with principal part $(X_0,0)$ over $S$.   The
notations of {\bf fine} and {\bf coarse moduli space} for the functor
Def$^-_{X_0}$ are defined in the same manner as for the functor Unf$^-_{f_0}$
in \S 1.   The objects we are going to classify are elements of

\begin{tabular}{lp{10cm}}
Def$^-_{X_0}(pt) =$ & $\{$ isomorphism classes of complex space germs $(X,0)$
which are isomorphic to a semiquasihomogeneous hypersurface singularity with
principal part $X_0\}$.
\end{tabular}

Again, as for Unf$^-_{f_0}$, we cannot expect to obtain fine moduli spaces in
general.
In order to obtain a coarse moduli space, we have to stratify $T_-$ into
$G$--invariant strata on which the geometric quotient with respect to $G$
exists, where $G = \exp\, L_+ \rtimes (E_{f_0} \cdot \C^\ast) \subset \mbox{
Aut}(T_-)$.  Once we have this, the proof is the same as for Theorem
\ref{1.3}.

We want to apply Theorem 4.7 from \cite{GP 2} to the action of $L_+$ on
$T_-$.

\begin{theorem}\label{4.1} (\cite{GP 2})\quad  Let $A$ be a noetherian
$\C$--algebra and $L_+ \subseteq \mbox{ Der}_\C^{nil}A$ a finite dimensional
nilpotent Lie algebra.   Suppose $A$ has a filtration

\[
F^\bullet :\; 0 = F^{-1}(A) \subset F^0 (A) \subset F^1 (A) \subset \ldots
\]

by subvector spaces $F^i(A)$ such that\\

$({\bf F})\qquad\qquad\qquad\qquad\qquad  \delta F^i(A) \subseteq F^{i-1} (A)\,
\mbox{ for all } i \in \Z,\; \delta \in L_+$.\\

Suppose, moreover,  $L_+$ has a filtration

\[
Z_\bullet : L_+ = Z_1 (L_+) \supseteq Z_2 (L_+) \supseteq \ldots \supseteq
Z_e(L_+) \supseteq Z_{e+1}(L_+) = 0
\]

by sub Lie algebras $Z_j(L_+)$ such that\\

$({\bf Z})\qquad\qquad\qquad\qquad\qquad\qquad
[L_+,  Z_j(L_+)] \subseteq Z_{j+1}(L_+)\, \mbox{for all}\, j \in
\Z$.\\

Let $d : A \to \mbox{ Hom}_\C (L_+, A)$ be the differential defined by $d(a)
(\delta) = \delta (a)$ and let \hbox{Spec $A = \cup U_\alpha$} be the
flattening
stratification of the modules

\[ \mbox{Hom}_\C (L_+, A) / A d (F^i(A))\quad i = 1, 2,\ldots \]

and

\[ \mbox{Hom}_\C (Z_j(L_+), A) / \pi_j(A(dA))\quad j = 1, \ldots, e, \]

where $\pi_j$ denotes the projection Hom$_\C(L_+, A) \to$
 Hom$_\C(Z_j(L_+), A)$.\\
Then $U_\alpha$ is invariant under the action of $L_+$ and $U_\alpha
\to U_\alpha / L_+$ is a geometric quotient which is a principal fibre bundle
with fibre $\exp(L_+)$.   Furthermore, the closure $\bar{U}_\alpha$ of
$U_\alpha$ is affine, $\bar{U}_\alpha = \mbox{ Spec } A_\alpha$, and the
canonical map $U_\alpha/L_+ \to \mbox{ Spec } A_\alpha^{L_+}$ is an open
embedding.
\end{theorem}

To apply the theorem we have to construct these filtrations and interpret the
corresponding stratification in terms of the Hilbert function of the Tjurina
algebra.\\

There are natural filtrations $H^\bullet (\C\{x\})$
respectively $F^\bullet (A_-)$ on $\C\{x\}$ respectively $A_-$ defined
as follows:\\

Let $F^i(A_-) \subseteq A_-$ be the $\C$--vectorspace
generated by all quasihomogeneous polynomials of degree $> - (i+1)w$ and
$H^i(\C\{x\})$ be the ideal generated by all quasihomogeneous polynomials
of degree $\ge i w$, where

\[w := \min\{w_1, \ldots, w_n\}.\]

The filtration $F^\bullet(A_-)$ has the property (${\bf F}$) because every
homogeneous vector field of $L_+$ is of degree $\ge w$.
We also have $A_- dA_- = A_- dF^sA_-$ with
$s = \left[\frac{(n-1)d-2\sum w_i}{w}\right]$, since $nd - 2 \sum w_i$
is the degree of the Hessian of $f$ and $t_k$ is the variable of smallest
degree.

To define $Z_\bullet$ let $Z_i(L_+) :=$ the Lie algebra generated by the
vectorfields $\delta \in L_+,\; \delta$ homogeneous and $\deg(\delta) \ge r_i$,
\[
r_i := \min\{\deg(\delta_j) \mid t_{k+1-j} \in F^{s-i}(A_-)\}.
\]
$Z_\bullet(L_+)$ has the property $({\bf Z})$ because
$\deg([\delta,\delta']) \ge
\deg(\delta) + \deg(\delta')$ for all $\delta, \delta' \in L_+$.

\begin{example}\label{4.2}
{\rm We continue with Example \ref{3.4}, $f_0 = x^3 + y^3 + z^7$.

$w = 3$.\\
$F^\circ(A_-)$ is the $\C$--vector space generated by $t_1, t_2, t_3, t_1^2,
t_1t_2, t_2^2$.\\
$F^1(A_-)$ is the $\C$--vector space generated by $t_4, \{t_1^\nu
t_2^\mu t_3^\lambda\}_{\nu+\mu+2\lambda \le 5}$.\\
$F^2(A_-)$ is the $\C$--vector space generated by $t_5, \{t_1^\nu
t_2^\mu t_3^\lambda t_4\}_{\nu+\mu+2\lambda \le 3}, \{t_1^\nu t_2^\mu
t_3^\lambda\}_{\nu + \mu + 2\lambda \le 8}$.\\
We have $s = 2 = \left[\frac{2 \cdot 21 - 2 \cdot 17}{3}\right]$.\\
$A_-dF^\circ(A_-) = \bigoplus\limits^3_{i=1} A_-dt_i$.\\
$A_-dF^1    (A_-) = \bigoplus\limits^4_{i=1} A_-dt_i$.\\
$A_-dF^2    (A_-) = A_-dA_-$.\\
$r_1 = 3, r_2 = 6$.\\
$L_+ = Z_1(L_+)$.\\
$Z_2(L_+)$ generated by the homogeneous vector fields $\delta \in L_+$ with
$\deg(\delta) \ge 6$.\\
Especially $A_-Z_2(L_+) = \sum\limits^5_{i=3} A_- \delta_i$.\\
$Z_3(L_+) = 0$.
}
\end{example}

We can use Theorem \ref{4.1} to obtain a geometric quotient of the action of
$L_+$ on the flattening stratification defined by the filtrations $F^\bullet$
and $Z_\bullet$.   Before doing this we shall prove that this flattening
stratification is also the flattening stratification of the modules defining
the Hilbert function of the Tjurina algebra.\\

For $t \in T_-$ the {\bf Hilbert function of the Tjurina algebra}
\[
\C\{x\}/\left(F(t),\frac{\partial F(t)}{\partial x_1}, \ldots, \frac{\partial
F(t)}{\partial x_n}\right)
\]
corresponding to the singularity defined by $t$ with respect to $H^\bullet$ is
by definition the function
\[
m \mapsto \tau_m(t) := \dim_\C \C\{x\}/\left(F(t), \frac{\partial
F(t)}{\partial x_1}, \ldots, \frac{\partial F(t)}{\partial x_n}, H^m\right).
\]

Notice that $\tau_m(t) = \tau(t)$ if $m$ is large and $\tau_m(t)$ does not
depend on $t$ for small $m$.   On the other hand, $\mu_m := \mu_m(t) :=
\dim_\C \C\{x\}/(\frac{\partial F(t)}{\partial x_1}, \ldots, \frac{\partial
F(t)}{\partial x_n}, H^m)$ does not depend on $t \in T_-$ and
\[
\mu_m - \tau_m(t) = \mbox{ rank } (\delta_i(t_j) (t))_{\deg(t_j) > d- mw}.
\]

This is an immediate consequence of the following fact:\\
Let
\[
T^m := A_-\{x\}/\left(F, \frac{\partial F}{\partial x_1}, \ldots,
\frac{\partial F}{\partial x_n}, H^m\right),
\]
then the following sequence is exact and splits:  let $\{X^\alpha\}_{\alpha
\in B}$ be a monomial base of $A_-\{x\}/\left(\frac{\partial F}{\partial x_1},
\ldots, \frac{\partial F}{\partial x_n}\right)$.
\[
\begin{array}{rcrllllll}
0 & \to & \bigoplus\limits_{\buildrel |\alpha|\le d\over \alpha \in B}
		 A_-x^\alpha & \to & T^{\frac{d}{w} + i} & \to
& \mbox{Der}_\C A_-/\left(\kl + \sum_{\deg(t_j) \le - iw} A_-
\frac{\partial}{\partial t_j}\right) & \to & 0\\
& &  x^\alpha & \mapsto & \mbox{class}(x^\alpha) & & & & \\
& & & &\mbox{class}(m_j) & \mapsto & \mbox{class}(\frac{\partial}{\partial
t_j}), & &
\end{array}
\]
and with the identification $\sum\limits_{\deg(t_j) > -iw} A_-
\frac{\partial}{\partial t_j} \simeq A_-^{N_i}$ we obtain\\
Der$_\C A_-/(\kl + \sum\limits_{\deg(t_j) \le -iw} A_- \frac{\partial}{\partial
t_j}) \simeq A_-^{N_i}/M_i$, where $M_i$ is the $A_-$--submodule
generated by the rows of the matrix $(\delta_i(t_j))_{\deg(t_j) >
-iw}$.

We have $F \in H^m$, hence $\mu_m = \tau_m$, if $m \le \frac{d}{w}$
and $H^m \subset
(\frac{\partial F}{\partial x_1}, \ldots, \frac{\partial F}{\partial x_n})$,
hence $\mu_m - \tau_m(t)$ is independent of $m$ and equal to $\mu - \tau(t)$,
if $m \ge \frac{d}{w} + s + 1$ .

Therefore, we have $s + 1$ relevant values for $\tau_i$, and we denote
\vspace{-0.5cm}

\begin{eqnarray*}
\underline{\tau}(t) & := & (\tau_{\frac{d}{w} + 1}(t), \ldots,
\tau_{\frac{d}{w} + s + 1}(t)),\\
\underline{\mu} & := & (\mu_{\frac{d}{w}+1}, \ldots, \mu_{\frac{d}{w}+s+1}).
\end{eqnarray*}

Moreover, let $\Sigma = \{\underline{r} := (r_1, \ldots, r_{s+1}) \mid
\exists\, t
\in T_-$ so that $\underline{\mu} - \underline{\tau}(t) = \underline{r}\}$ and
$T_- = \cup_{\underline{r} \in \Sigma}U_{\underline{r}}$ be the flattening
stratification of the modules
$T^{\frac{d}{w}+1}, \ldots, T^{\frac{d}{w}+s+1}$.   That is,
$\{U_{\underline{r}}\}$ is the stratification of $T_-$ defined by fixing the
Hilbert function $\underline{\tau} = \underline{\mu} - \underline{r}$ with the
scheme structure defined by the flattening property.

Let us now consider an arbitrary deformation $\phi : (\kx,\{0\} \times S)
\hookrightarrow (\C^n,0) \times S \to S$ of $(X,0) \subset (\C^n,0)$ of
negative weight over a base space $S \in \ks$ where, for each $s \in S$, the
ideal of the germ $(\kx,(0,s)) \subset (\C^n \times S, (0,s))$ is defined by
$F(x,s) = f(x) + g(x,s),\; g(x,0) = 0$.

Let us denote by $\ko_S \{x\} = \ko_{\C^n \times S, 0 \times S}$ the
topological restriction of $\ko_{\C^n \times S}$ to $0 \times S$,
considered as a sheaf on $S$.   Then
$J(I_{\kx , 0 \times S})$, the Jacobian ideal sheaf of $(\kx, \{0\} \times S)
\subset (\C^n,0) \times S$, is locally defined by $(F, \frac{\partial
F}{\partial x_1}, \ldots, \frac{\partial F}{\partial x_n}) \subset \ko_S
\{x\}$ and $H^m_S \subset \ko_S\{x\}$ is the ideal sheaf generated by $g \in
\ko_S \{x\}$ such that $\deg_x g \ge mw,\; w = \min\{w_1, \ldots, w_n\}$ as
above.   We say that the {\bf family} $\phi$ is {\bf
$\underline{\tau}$--constant} if the coherent $\ko_S$--sheaves
\[
T^m_S := \ko_S \{x\}/J(I_{\kx,\{0\}\times S}) + H^m_S
\]
are flat for $\frac{d}{w} + 1\le m \le \frac{d}{w} + s + 1$ (equivalently, for
all $m$).   Of course, if $T^m_S$ is flat, then
\[
\tau_m(s) := \dim_\C T^m_{S,s} \otimes \ko_{S,s}/\frak m_{S,s}
\]
is independent of $s \in S$.   The converse holds for reduced base spaces:

\begin{lemma}
If $S$ is reduced, then the sheaf $T^m_S$ is flat if and only if $\tau_m(s)$ is
independent of $s \in S$.
\end{lemma}

The proof is standard (cf.\ \cite{GP 3}).   Hence, over a reduced base space
$S$, $\underline{\tau}$--constant means just that the Hilbert function
$\underline{\tau}(s) = (\tau_{\frac{d}{w} + 1} (s), \ldots, \tau_{\frac{d}{w}
+ s + 1} (s))$ of the Tjurina algebra is constant.   But for arbitrary base
spaces we have to require flatness of the corresponding $T^m_S$.

{\bf Example} ($f_0 = x^3 + y^3 + z^7$, continued)
\begin{eqnarray*}
\underline{\tau}(t) & = & (\tau_8(t), \tau_9(t), \tau_{10}(t))\\
\underline{\mu} & = & ( \mu_8, \mu_9, \mu_{10}) = (22, 23, 24)\\
\Sigma & = & \{(0, 0, 0), (0, 0,1), (0, 1, 2), (1, 1, 2), (1, 2, 3)\}\\
U_{(1,2,3)} & = & D(2t_3 - \frac{10}{7} t_1 t_2) \subseteq T_- = \C^5\\
U_{(1, 1, 2)} & = & V(2t_3 - \frac{10}{7}t_1t_2) \cap D(t_1, t_2) \subseteq
T_-\\
U_{(0, 1, 2)} & = & V(t_1, t_2, t_3) \cap D(t_4) \subseteq T_-\\
U_{(0, 0, 1)} & = & V(t_1, t_2, t_3, t_4) \cap D(t_5) \subseteq T_-\\
U_{(0, 0, 0)} & = & \{(0, 0, 0, 0, 0)\}.
\end{eqnarray*}

\begin{lemma}\label{4.3}
\begin{enumerate}
\item $(0, \ldots, 0,1)$ and $(0, \ldots, 0) \in \Sigma$.   $U_{(0, \ldots,0)}
=
\{0\}$ is a smooth point and $U_{(0, \ldots, 1)}$ is defined by $t_1 = \cdots
= t_{k-1} = 0$ and $t_k \not= 0$.

\item Let $\bar{\Sigma} = \Sigma\backslash\{(0, \ldots, 0)\}$ and for
$\underline{r} \in
\bar{\Sigma}$ put
\[
\widetilde{U}_{\underline{r}} = \left\{
\begin{array}{ll}
U_{\underline{r}} & \mbox{ if } \underline{r} \not= (0, \ldots, 0, 1)\\
U_{(0, \ldots, 0, 1)} \cup U_{(0, \ldots, 0)} & \mbox{ if } \underline{r} =
(0, \ldots, 0, 1).
\end{array}\right.
\]
\end{enumerate}

Then $\{\widetilde{U}_{\underline{r}}\}_{\underline{r} \in \bar{\Sigma}}$ is
the
flattening stratification of the modules $\{\mbox{Hom}_\C(L_+, A_-)/A_- dF^i
A_-\}$ and $\{\mbox{Hom}_\C(Z_i(L_+), A_-)/\pi_i(A_-dA_-)\}$.
\end{lemma}

{\bf Proof of Lemma \ref{4.3}}:  Because of the exact sequence above the
flattening stratification of the modules $\{T^{\frac{d}{w}+i}\}$ is also the
flattening stratification of $\{\mbox{Der}_\C A_-/(\kl + \sum_{\deg(t_j) \le
-iw} A_- \frac{\partial}{\partial t_j})\}$ respectively the flattening
stratification of $\{A_-^{N_i}/M_i\},\; M_i$ the submodule generated by the
rows of the matrix $(\delta_i(t_j))_{\deg(t_j) > -iw}$.

Now we have

$(\ast)$\hspace{4.5cm}$\delta_i(t_j) = \delta_{k-j+1}(t_{k-i+1})$.

By definition of $Z_i(L_+)$ we have
\[
A_-Z_i(L_+) = \sum_{t_{k+1-j} \in F^{s-i}} A_- \delta_j
\]
and with the identification
\[
\sum\limits A_- \frac{\partial}{\partial t_j} = A_-^k,
\]
and $M^i$ the submodule generated by the rows of the matrix
$(\delta_\ell(t_j))_{\ell \ge r_i}$ we obtain
\[
\mbox{Der}_\C A_- / A_- Z_i(L_+) \cong A_-^k/M^i.
\]

(*) implies that the flattening stratification of the modules
$\{T^{\frac{d}{w}+1}, \ldots, T^{\frac{d}{w}+s}\}$, which is $T_- =
\cup_{\underline{r} \in
\bar{\Sigma}} \widetilde{U}_{\underline{r}}$, is the flattening stratification
of the
modules $\{\mbox{Der}_\C A_-/A_- Z_i(L_+)\}_{i=1, \ldots, s}$.

Furthermore the modules $\{\mbox{Hom}_\C(L_+, A_-)/A_-dF^iA_-\}$ and\\
$\{\mbox{Der}_\C A_-/A_- L_+ + \sum_{\deg(t_j)\le -iw}
A_-\frac{\partial}{\partial t_j}\}$ have the same flattening stratification
and they are flat on $U_{\underline{r}}$, because
\[
0 \to A_- \to \mbox{ Der}_\C A_-/A_-L_+ + \sum_{\deg(t_j)\le -iw}
A_-\frac{\partial}{\partial t_j}   \to \mbox{ Der}_\C A_- /\kl +
\sum_{\deg(t_j)\le -iw} A_-\frac{\partial}{\partial t_j} \to 0
\]
is exact and splits on $T_-\backslash\{0\}$.

This proves the lemma.

\begin{remark}\label{4.4}{\rm
The main point of the lemma is that the flattening stratification of the
modules $\{\mbox{Hom}_\C(L_+, A_-)/A_- dF^iA_-\}$ is equal to the
flattening stratification of the modules $\{\mbox{Hom}_\C(Z_i(L_+),
A_-)/\pi_i(A_- d A_-)\}$, hence, is defined by the Hilbert function of the
Tjurina algebra alone, without any reference to the action of $L$.   This is a
consequence of the symmetry expressed in Proposition \ref{3.2}.}
\end{remark}

As a corollary we obtain the following

\begin{theorem}\label{4.5}
For $\underline{r}
\in \Sigma,\; \widetilde{U}_{\underline{r}}$ is invariant under the action of
$L_+$.
Let Spec $A_{\underline{r}}$ be the closure of $\widetilde{U}_{\underline{r}}$
then
$\widetilde U_{\underline{r}} \to \widetilde{U}_{\underline{r}}/L_+$ is a
geometric
quotient contained in Spec $A^{L_+}_{\underline{r}}$ as an open subscheme
of Spec $A^{L_+}_{\underline{r}}$.
\end{theorem}

{\bf Example} ($f_0 = x^3 + y^3 + z^7$, continued)
\[
\begin{array}{lccc}
1) & \widetilde{U}_{(1, 2, 3)} = D(2t_3 - \frac{10}{7} t_1 t_2) &
\longrightarrow &
\widetilde{U}_{(1, 2, 3)}/L_+ = \mbox{ Spec}\, \C[t_1, t_2,
t_3]_{2t_3-\frac{10}{7} t_1t_2}\\[1.0ex]
& \bigcap\mid & & \bigcap\mid\\[1.0ex]
& \mbox{Spec }\C[t_1, \ldots, t_5] & \longrightarrow & \mbox{Spec }\C[t_1,
t_2, t_3]\\[2.0ex]
2) & \widetilde{U}_{(1, 1, 2)} & \longrightarrow & \widetilde{U}_{(1, 1,
2)}/L_+ =
D(t_1, t_2)\\[1.0ex]
& \bigcap\mid & & \bigcap\mid\\[1.0ex]
& \mbox{Spec }\C[t_1, t_2, t_4, t_5] & \longrightarrow & \mbox{Spec }\C[t_1,
t_2, t_4]\\[1.0ex]
\multicolumn{4}{l}{(\mbox{identifiying } \C[t_1, \ldots, t_5]/2t_3 -
\frac{10}{7} t_1t_2 = \C[t_1, t_2, t_4, t_5].)}\\[2.0ex]
3) & \widetilde{U}_{(0,1,2)} & \longrightarrow & \widetilde{U}_{(0,1,2)}/L_+ =
D(t_4)\\[1.0ex]
& \bigcap\mid & & \bigcap\mid\\[1.0ex]
& \mbox{Spec }\C[t_4, t_5] & \longrightarrow & \mbox{Spec }\C[t_4]\\[2.0ex]
4) & \widetilde{U}_{(0,0,1)} & = & \widetilde{U}_{(0,0,1)}/L_+\\[1.0ex]
& \| & & \|\\[1.0ex]
& \mbox{Spec }\C[t_5] & = & \mbox{Spec }\C[t_5]
\end{array}
\]
\vspace{1cm}

Now $L/L_+ \simeq \C \delta_1$ acts on the geometric quotients
$\widetilde{U}_{\underline{r}}/L_+$ (the $\C^\ast$--action defined by the Euler
vector field $\delta_1$).   Also the group $E_{f_0}$ acts and this action
commutes
with the $\C^\ast$--action (cf.\ \ref{2.5}).   If we combine this fact with
Theorem 4.6 we obtain the main theorem of this article.   In order to
formulate it properly let us denote by
\[
\mbox{Def}^-_{X_0, \underline{\tau}} : \ks \to \mbox{ sets}
\]
the subfunctor of Def$^-_{X_0}$ which associates to a base space $S \in \ks$
the set of isomorphism classes of $\underline{\tau}$--constant families of
deformations of negative weight with principal part $(X_0,0)$ over $S$.   For
such a family $\underline{\tau}(s)$ is constant and equal to some tuple
$\underline{\mu} - \underline{r} \in \N^{s+1}$.

\begin{theorem}\label{4.6}
Let $G = \exp L_+ \rtimes (E_{f_0} \cdot \C^\ast) \subseteq \mbox{ Aut }(T_-)$.

\begin{enumerate}
\item The orbits of $G$ are unions of finitely many integral manifolds of
$\kl$.

\item Let $T_- = \cup_{\underline{r}\in \Sigma} U_{\underline{r}}$ be the
stratification fixing the
Hilbert function $\underline{\tau}$ of the Tjurina algebra described above.
$U_{\underline{r}}$ is invariant under the action of $G$ and the geometric
quotient $U_{\underline{r}} \to U_{\underline{r}}/G$ exists
and is locally closed in a weighted projective space.

\item $U_{\underline{r}}/G$ is the coarse moduli space for the functor
Def$^-_{X_0, \underline{\tau}} : $ complex spaces $\to$ sets with
$\underline{\tau} = \underline{\mu} - \underline{r}$.
\end{enumerate}
\end{theorem}

\begin{remark}{\rm
As in the case of right equivalence (see Remark 1.5) we may take (separated)
algebraic spaces as category of base spaces.   That is, $U_{\underline{r}}/G$
is a coarse moduli space for the functor
\[
\mbox{Def}^-_{X_0, \underline{\tau}} : \mbox{ algebraic spaces }\to \mbox{
sets}.
\]
}
\end{remark}

{\bf Proof} (of Theorem 4.7):  We first prove that $U_{\underline{r}}$ is
invariant under the
action of $G$ and that $U_{\underline{r}} \to U_{\underline{r}}/G$ is a
geometric quotient.

To prove that $U_{\underline{r}}$ is invariant under the action of $G$ it is
enough by definition of $U_{\underline{r}}$ that it is invariant under the
action of $E_{f_0}$.   The Hilbert function $\underline{\tau}$ of the Tjurina
algebra is invariant under contact equivalence.   This is a consequence of
Theorem \ref{2.1} because an automorphism $\varphi$ of $\C\{x\}$ inducing the
isomorphy of two semiquasihomogeneous singularities with principal part $f_0$
has degree $\ge 0$.   More precisely, let $f,g$ be semiquasihomogeneous with
principal part $f_0$ and $uf = \varphi(g)$ for a unit $u$ then $\deg(\varphi)
\ge 0$ and consequently $(f, \frac{\partial f}{\partial x_1}, \ldots,
\frac{\partial f}{\partial x_n}, H^m)$ is mapped isomorphically to
$(g, \frac{\partial g} {\partial x_1}, \ldots, \frac{\partial g}{\partial x_n},
H^m)$ for all $m$, in particular $\underline{\tau}(f) = \underline{\tau}(g)$.

Moreover, let $\sigma \in E_{f_0}$, then there is a $\varphi : A_-\{x\} \to
A_-\{x\},\; \deg_x(\varphi) \ge 0$ and $\varphi|_{A_-} = id_{A_-}$
such that
\[\varphi(F(x,t)) \equiv F(x,\sigma (t))\hbox{ mod } A_- H^N
 \hbox{ \ for sufficiently large }N\]
(cf. proof of Proposition 2.4).

This implies $\sigma(T^m) = T^m$ for all $m$ and proves that $E_{f_0}$ and,
therefore, $G$ acts on the strata $U_{\underline{r}}$ of the flattening
stratification of the modules $\{T^m\}$.

Now we prove that $U_{\underline{r}} \to U_{\underline{r}}/G$ is a geometric
quotient.   First of all it is obvious that the geometric quotients
\[
U_{(0, \ldots, 0,1)} \to U_{(0, \ldots, 0,1)}/G = \{\ast\}
\]
and
\[
U_{(0, \ldots, 0)} = \{\ast\} = U_{(0, \ldots, 0)}/G = \{\ast\}
\]
exist.

Let $\underline{r} \not= (0, \ldots, 0,1),\; (0,\ldots, 0)$ then
$\widetilde{U}_{\underline{r}} = U_{\underline{r}}$.   Let
$U_{\le\underline{r}} =
\mbox{ Spec}A_{\underline{r}}$ be the closure of $U_{\underline{r}}$ then we
obtain
\[
\begin{array}{ccc}
\mbox{Spec}A_{\underline{r}} & \buildrel \pi\over\longrightarrow &
\mbox{Spec}A_{\underline{r}}^{L_+}\\[1.0ex]
\cup|\; i & & \cup|\; j\\[1.0ex]
U_{\underline{r}} & \buildrel\pi|_{U_{\underline{r}}}\over\longrightarrow &
U_{\underline{r}}/L_+.
\end{array}
\]

$\pi|_{U_{\underline{r}}}$ defines a geometric quotient and $i,j$ are open
embeddings (Theorem \ref{4.5}).   Notice that $\pi$ itself is not necessarily a
geometric quotient.

Now Spec$A_{\underline{r}}^{L_+}$ is affine and $E_{f_0}$ acts on
Spec$A_{\underline{r}}^{L_+}$ and also on $U_{\underline{r}}/L_+$.   This
implies (cf.\ \cite{MF}) that
\[
\mbox{Spec}A_{\underline{r}}^{L_+} \buildrel\lambda\over\to \mbox{
Spec}(A_{\underline{r}}^{L_+})^{E_{f_0}}
\]
is a geometric quotient (not necessarily as algebraic schemes since
$A^{L_+}_{\underline{r}}$ need not be of finite type over $\C$) and
consequently
\[
\lambda|_{U_{\underline{r}}/L_+} : U_{\underline{r}}/L_+ \to
(U_{\underline{r}}/L_+)/E_{f_0}
\]
is a geometric quotient which is an algebraic scheme.   Especially
$(U_{\underline{r}}/L_+)/E_{f_0} \subseteq
\mbox{ Spec}(A_{\underline{r}}^{L_+})^{E_{f_0}}$ is an open subset.

Finally, $\C^\ast$ acts on Spec$(A_{\underline{r}}^{L_+})^{E_{f_0}}$.   It has
one
fixed point $\{\ast\}$ corresponding to $U_{(0, \ldots, 0)} \subseteq
\bar{U}_r =$ Spec$A_{\underline{r}}$.   Outside this fixed point the
$\C^\ast$--action leads to a geometric quotient:
\[
\begin{array}{ccc}
\mbox{Spec}(A_{\underline{r}}^{L_+})^{E_{f_0}}\backslash\{\ast\} &
\longrightarrow
& \mbox{Proj}(A_{\underline{r}}^{L_+})^{E_{f_0}}\\[1.0ex]
\cup & & \cup\\[1.0ex]
(U_{\underline{r}}/L_+)/E_{f_0} & \longrightarrow &
((U_{\underline{r}}/L_+)/E_{f_0})/\C^\ast\\[1.0ex]
& & \|\\[1.0ex]
& & U_{\underline{r}}/G.
\end{array}
\]
This proves part (1) and (2) of the theorem.

It remains to prove that if $t,t' \in T_-$
define isomorphic singularities then
$t$ and $t'$ are in the same orbit of $G$.

Let $F_t = u\varphi(F_{t'})$ for $t, t' \in T_-, u \in \C\{x\}^\ast$ a unit and
$\varphi$ an automorphism of $\C\{x\}$.
Using the $\C^\ast$--action we find
$t''\in T_-,\ u_1 =\frac{u}{u(0)}\in \C\{x\}^\ast$ and an automorphism
$\varphi_1$ of $\C\{x\}$ such that $F_t = u_1\varphi_1 (F_{t''}),\ u_1(0)=1$
and $t'$ and $t''$ are in one $\C^\ast$--orbit. Then
\[
G(z) := (1+z(u_1-1))\varphi_1(F_{t''})
\]
is an unfolding of $G(0) = F_t$ of negative weight.
This unfolding can be induced
by the semiuniversal unfolding, that is there exists a family of coordinate
transformations $\underline{\psi}(z, -)$ and a path $v$ in $T_-$ such that
\[
G(z) = F(\psi_1(z,x), \ldots, \psi_n(z,x), v(z))
\]
and $v(0) = t$ and $F_{t''} \buildrel r\over\sim F(\psi(1,x), v(1))$.
Now $t=v(0)$ and $v(1)$ are in one orbit of $\exp L$, and $v(1)$ and $t''$
are in one orbit of $E_{f_0}$.
Hence the result.

Now (3) follows in the same manner as the proof of Theorem 1.3.

{\bf Example} ($f_0 = x^3 + y^3 + z^7$, continued)

\begin{enumerate}
\item $U_{(1,2,3)} \longrightarrow U_{(1,2,3)}/G \simeq \C^2,\;
\underline{\tau} = (21, 21, 21),\; \tau = 21$ \\
normal form:  $f_0 + t_1 xz^5 + t_2 yz^5 + t_3xyz^3,\\
(t_1:t_2:t_3) \in D_+ (2t_3 - \frac{10}{7} t_1t_2)/S_3 \subset
\P^2_{(1:1:2)}/S_3$\\
$(D_+(2t_3 - \frac{10}{7} t_1t_2)/S_3 \simeq \C^2$, the $S_3$--action being
explained in Example 2.8).

\item $U_{(1,1,2)} \longrightarrow U_{(1,1,2)}/G \simeq
\P^2_{(2,3,5)}\backslash(0:0:1),\; \underline{\tau} = (21, 22, 22),\; \tau =
22$\\
normal form:  $f_0 + t_1xz^5+t_2yz^5+\frac{10}{7} t_1t_2xyz^3 +t_4 xyz^4$,\\
$(t_1:t_2:t_4) \in \P^2_{(1:1:5)}/S_3 \;\; (\simeq \P^2_{(2,3,5)})$

\item $U_{(0,1,2)} \longrightarrow U_{(0,1,2,)}/G = \{\ast\},\;
\underline{\tau} = (22, 22, 22),\; \tau = 22$\\
normal form: $f_0 + xyz^4$

\item $U_{(0,0,1)} \longrightarrow U_{(0,0,1,)}/G = \{\ast\},\;
\underline{\tau} = (22, 23, 23),\, \tau = 23$\\
normal form:  $f_0 + xyz^5$

\item $U_{(0,0,0)} \longrightarrow U_{(0,0,0,)}/G = \{\ast\},\;
\underline{\tau} = (22, 23, 24),\; \tau = 24$\\
normal form:  $f_0$
\end{enumerate}

Hence the moduli space of semiquasihomogeneous hypersurface singularities $X =
\{(x,y,z) \mid f (x,y,z) = 0\}$ with principal part $X_0 = \{(x,y,z) \mid x^3
+ y^3 + z^7 = 0\}$ consists of 5 strata ($\C^2,\; \P^2_{(2,3,5)} \backslash (0
: 0 : 1)$, and 3 isolated points) corresponding to 5 possible Hilbert
functions $\underline{\tau}$ of the Tjurina algebra $\C\{x,y,z\}/(f,
\frac{\partial f}{\partial x}, \frac{\partial f}{\partial y}, \frac{\partial
f}{\partial z})$.   The generic stratum $U_{(1,2,3)}$ (minimal
$\underline{\tau}$) is an open subset in $\C^5$, the quotient being
2--dimensional, as well as the quotient of the 4--dimensional ``subgeneric''
stratum $U_{(1,1,2)}$.
Note that the families of normal forms are not universal.   It just means that
each semiquasihomogeneous singularity with principal part $f_0$ occurs and
that different parameters do not give contact equivalent singularities, except
modulo the $\C^\ast$-- and $S_3$--action.

We see that $U_{(1,1,2)}/G$ can be compactified by $U_{(0,1,2)}/G$, that is
\[
U_{(1,1,2)} \cup U_{(0,1,2)} \to (U_{(1,1,2)} \cup U_{(0,1,2)})/G =
\P^2_{(2,3,5)}
\]
is a geometric quotient. So in this example there exist
geometric quotients of the strata with constant Tjurina number and, hence, a
coarse moduli space for fixed principal part and fixed Tjurina number.
In general this is false (cf.\ \cite{LP}, \S 7).

\begin{remark}{\rm
1.\quad The generic stratum $U_{\underline{\tau}\min}$ corresponding to minimal
Hilbert function $\underline{\tau}$ (with respect to lexicographical ordering)
is an open, quasiaffine subset of $T_-$ and, hence,
$U_{\underline{\tau}\min}/L_+$ is smooth by Theorem 4.1.   In particular, the
generic moduli space $U_{\underline{\tau}\min}/G$ has, at most, quotient
singularities (coming from the $\C^\ast$--action and the finite group
$E_{f_0}$).
It is not known whether the bigger stratum $U_{\tau\min}$ corresponding to
minimal Tjurina number $\tau$ admits a geometric quotient, except for $n = 2$
(cf.\ \cite{LP}).

2.\quad We always have two special strata, the most special $U_{(0, \ldots,0)}
=
\{\ast\}$ (corresponding to $f_0$) and the ``subspecial'' $U_{(0, \ldots, 1)}
\cong \C\backslash \{\ast\}$ (corresponding to the singularity $f_0 + m_k,\;\;
m_k$ generating the socle of $\C\{x\}/j(f_0)$, that is the monomial of
maximal degree).   The $G$--quotients of these strata give two reduced,
isolated points.

3.\quad As we have seen for $x^3 + y^3 + z^7$, the finite group $E_{f_0}$ need
not
be abelian.   If $f_0 = x_1^{a_1} + \cdots + x_n^{a_n}$ is of Brieskorn--Pham
type and gcd$(a_i, a_j) = 1$ for $i \not= j$, then $E_{f_0} \cong \mu_d$, the
group of $d$'th roots of unity, $d = \deg\, f_0$.

4.\quad Note that a coarse moduli space is more than just a bijection between
its points and the corresponding set of isomorphism classes.   For instance,
let $U_{\underline{r}}/G$ be affine and let $S \buildrel\sigma\over\to (\kx,
\sigma(S)) \buildrel\phi\over\to S$ be a family of deformations from
Def$^-_{X_0} (S)$ with $\underline{\tau}(\kx_s, \sigma(s)) = \underline{\mu} -
\underline{r}$.   If $S$ is compact then $\phi$ must be locally trivial since
any morphism from $S$ to $U_{\underline{r}}/G$ maps $S$ onto finitely many
points.}
\end{remark}

\end{document}